\documentclass[a4 paper,11pt]{article}

\usepackage[latin1]{inputenc}  
\usepackage[T1]{fontenc}  
\usepackage[english,frenchb]{babel}  
\usepackage{textcomp}  
\usepackage[all]{xy}  

\usepackage{amssymb,amsmath,latexsym,epsfig,euscript}

\def\codim {{\rm \mbox{codim}}}

\def\Conv {{ \mbox{\rm Conv}}}
\def\conv {{ \mbox{\rm  \scriptsize Conv}}}
\def\Vol {{\mbox{\rm Vol}}}

\def\ass {{\rm \mbox{Ass}}}

\def\Supp {{\mbox{\rm Supp\,}}}

\def\Res {\mbox{\rm R{\'e}s}}

\def\mmax {{\mbox{\rm \scriptsize sup}}}

\def\ess {{ \mbox{\rm \scriptsize ess}}}
\def\abs {{ \mbox{\rm \scriptsize abs}}}

\def\Div {\mbox{\rm Div}}
\def\div {\mbox{\rm div}}
\def\Tors {\mbox{\rm Tors}}

\def\lg {\mbox{\rm lg}}
\def\Ass {\mbox{\rm \scriptsize Ass}}
\def\Irr {\mbox{\rm Irr}}
\def\irr {\mbox{\rm \scriptsize Irr}}

\def\W {{\mbox{\rm \scriptsize W}}}

\def\NP {{\mathcal{N}}}
\def\Ch {{\mathcal{C}{\it h}}}
\def\hnorm {{\widehat{h}}}
\def\arith {{{\mathcal A}, \alpha}}

\def\Card {\mbox{\rm Card}}
\def\Faces {F} 

\def\A{\mathbb{A}}
\def\C{\mathbb{C}}
\def\G{\mathbb{G}}
\def \N{\mathbb{N}}
\def \P{\mathbb{P}}
\def \Q{\mathbb{Q}}
\def \R{\mathbb{R}}
\def \T{\mathbb{T}}
\def \Z{\mathbb{Z}}

\def\cA {{\cal A}}
\def\cB {{\cal B}}

\def\cO {{\cal O}}

\def\cV {{\cal V}}

\def\ga {{\mathfrak{a}}}
\def\gf {{\mathfrak{f}}}
\def\gothg {{\mathfrak{g}}}
\def\gp {{\mathfrak{p}}}
\def\gq {{\mathfrak{q}}}

\def\Qbarra {{\overline{\Q}}}
\def\ov#1{{\overline{#1}}}

\newtheorem{lem}{Lemme}[section]

\newtheorem{prop}[lem]{Proposition}
\newtheorem{thm}[lem]{Th{\'e}or{\`e}me}
\newtheorem{cor}[lem]{Corollaire}

\newtheorem{defn}{D{\'e}finition}[section]

\newtheorem{rem}[defn]{Remarque}
\newtheorem{exmpl}[defn]{Example}

           {\vspace{3.3mm}
           \noindent{\bf #1}\it}%
           {\vspace{3.3mm}}

\newenvironment{demo}{\trivlist\item[\hskip 
\labelsep{\it D{\'e}monstration.--}]}{\hfill\mbox{$\square$} \endtrivlist}

\newenvironment{Demo}[1]{\trivlist\item[\hskip 
\labelsep{\it #1}]}{\hfill\mbox{$\square$} \endtrivlist}

\begin{document}

\noindent{\LARGE{\bf Minimums successifs des vari{\'e}t{\'e}s toriques 
projectives\footnote{Version du 2 Juin 2004.} 
}

\vspace{2mm}

\noindent {\large Mart{\'\i}n Sombra\footnote{Financ{\'e} par une 
bourse  
post-doctorale Marie Curie  
 du programme europ{\'e}en 
{\em Improving  Human
Research Potential and the Socio-economic Knowledge Base}, 
contrat n\textordmasculine \ HPMFCT-2000-00709. 
}}
}

\vspace{7mm}

% ----------- Abstract -----------------------------------------------

\noindent {\small {\bf R{\'e}sum{\'e}.--- \ }} 
{\small 
Nous calculons les minimums successifs de la
vari{\'e}t{\'e} torique 
projective 
$X_\cA$ associ{\'e}e {\`a} un ensemble fini $ \cA \subset \Z^n$. 
Comme cons{\'e}quence de ce calcul et 
des r{\'e}sultats  de S.-W. Zhang sur la r{\'e}partition des petits
points, nous d{\'e}duisons des estimations pour la hauteur
de la sous-vari{\'e}t{\'e} $X_\cA$ et du
 $\cA$-r{\'e}sultant. 
Ces estimations nous  permettent d'obtenir un analogue arithm{\'e}tique
du th{\'e}or{\`e}me de B{\'e}zout-Koushnirenko
sur le nombre de solutions d'un syst{\`e}me d'{\'e}quations
polynomiales.  

Comme application de ce r{\'e}sultat, nous am{\'e}liorons
les estimations connues pour la hauteur des polyn{\^o}mes dans le 
Null\-ste\-llen\-satz creux. 

}

\medskip

\noindent {\small {\bf Abstract.--- \ \ Successive minima of projective 
toric varieties.}} 
{\small
We compute the successive minima of the projective toric 
variety 
$X_\cA$ associated to a finite set 
$ \cA \subset \Z^n$. 
As a consequence of this computation and of the results of
S.-W. Zhang on the distribution of small points, 
we derive estimates
for the height of the subvariety
$X_\cA$ and of the 
 $\cA$-resultant. 
These estimates allow us to obtain an arithmetic analogue of 
the B{\'e}zout-Kushnirenko's theorem concerning the number of solutions of a
system of polynomial equations.
 
As an application of this result, we improve the known estimates
for the height of the polynomials in the sparse Nullstellensatz. 

}

\medskip 

\noindent {\small {\bf Mots clefs.--- \ } Hauteur,
minimums successifs, vari{\'e}t{\'e} torique, r{\'e}sultant creux, 
th{\'e}or{\`e}me de Koush\-ni\-renko, Nullstellensatz arithm{\'e}tique.}
 
\medskip

\noindent {\small {\bf Classification math{\'e}matique par sujets
    (2000).--- \ } 
{Primaire:} 11G50; 
{Secondaire:} 14G40, 14M25.}

% ----------- Intro -----------------------------------------------

\typeout{Intro}

\section*{Introduction et r{\'e}sultats}

L'{\'e}tude de la  r{\'e}partition 
des points de petite hauteur (ou petits points)
d'une vari{\'e}t{\'e} alg{\'e}brique a
re{\c c}u une attention consid{\'e}rable  
au cours des
derni{\`e}res ann{\'e}es, autour du probl{\`e}me de Bogomolov 
et de  ses 
g{\'e}n{\'e}ralisations~\cite{Ullmo98}, \cite{Zhang95}, 
\cite{BoZa95}, \cite{Bilu97}, \cite{DaPh98},
 \cite{DaPh99}, 
\cite{DaPh01}, 
{\it voir} {\'e}galement~\cite{Dav03} pour un aper{\c c}u historique. 
La notion de minimums successifs a {\'e}t{\'e}
introduite dans ce contexte par Zhang, qui  a aussi montr{\'e} leur  
{\'e}troite relation 
avec la hauteur de la vari{\'e}t{\'e} en question. 

\bigskip

Ici nous calculons 
les minimums successifs d'une vari{\'e}t{\'e} torique projective, d'o{\`u} nous
d{\'e}duisons 
des estimations pour sa hauteur comme sous-vari{\'e}t{\'e} et pour la hauteur 
du r{\'e}sultant creux.

Soit  \,$h$\,  la hauteur projective standard des points de 
\,$\P^N(\Qbarra)$,  
d{\'e}finie en utilisant la m{\'e}trique 
euclidi{e}nne  pour les places {\`a} l'infini. 
Soit  $V \subset \P^N$  une vari{\'e}t{\'e} quasi-projective, 
de dimension $r$.
Le {\em $i$-{\`e}me minimum successif} de $V$ (par rapport {\`a} la
hauteur projective $h$) est 
$$
\mu_i(V) := \sup \,
\left\{
\ \inf
\, \left\{ h(\xi) \ ; \ \xi \in V \setminus W
\right\} \ ; \ W \subset V, \ \codim_V (W ) = i
\right\}
$$
pour $i=1, \dots, r+1$, o{\`u} le supr{\^e}mum est pris sur toutes les
sous-vari{\'e}t{\'e}s $W$ de $V$ de codimension $i$. 
On {\'e}crit \,$
\mu^\ess(V):= \mu_{1}(V) $ \ 
et \ $ \mu^\abs(V):= \mu_{r+1} (V)
$\, 
pour les minimums {\em essentiel}  et
{\em absolu}, respectivement; 
on a 
$$
\mu^\ess(V) = \mu_{1}(V) \ge \cdots \ge \mu_{r+1}(V) = \mu^\abs(V).
$$

\smallskip

Soit \,$\cA =\{ a_0, \dots, a_N \} \subset \Z^n$\, un ensemble fini
de vecteurs entiers, et consid{\'e}rons l'application monomiale
$$
\varphi_\cA : \T^n \to \P^N
\quad \quad , \quad \quad
t \mapsto (t^{a_0} : \cdots : t^{a_N}). 
$$
o{\`u} 
\,$
\T^n:= (\Qbarra^\times)^n
$\, d{\'e}signe le tore alg{\'e}brique sur $\Qbarra$ de dimension $n$. 
La  
{\em
vari{\'e}t{\'e} torique projective}  \,$X_\cA \subset \P^N$\,  
associ{\'e}e est d{\'e}finie 
comme {\'e}tant la cl{\^o}ture de Zariski de
l'ensemble image de cette application
$$
X_\cA^\circ:= \varphi_\cA(\T^n) = \Big\{   \,(t^{a_0} : \cdots : t^{a_N})
\ ; \ t \in  \T^n \,\Big\} \ \subset \P^N. 
$$
Soit $L_\cA \subset \Z^n $ le sous-module engendr{\'e} par 
les diff{\'e}rences des vecteurs $a_0, \dots, a_N$; 
ceci est donc un {\it r{\'e}seau}
de l'espace lin{\'e}aire engendr{\'e} $L_\cA \otimes_\Z \R \subset \R^n$.  
On consid{\`e}re la forme volume $\Vol_\cA$ sur 
cet espace lin{\'e}aire, 
invariante par translations et telle que
$\Vol (S)= 1$ pour un  simplex {\'e}l{\'e}mentaire quelconque $S$  
de $L_\cA$.

Soit \,$Q_\cA \subset \R^n$ le polytope associ{\'e}, 
d{\'e}fini comme l'enveloppe convexe $\Conv(\cA)$ de l'ensemble $\cA$. 
Posons $\dim(\cA):= \dim (L_\cA)$ et $\Vol(\cA):=  \Vol_\cA(Q_\cA)$. 
Avec ces notations, la dimension et 
le degr{\'e} de $X_\cA $ s'explicitent comme la dimension et le
volume de l'ensemble 
$\cA$, respectivement: 
$$
\dim (X_\cA) = \dim(\cA)= \dim (L_\cA)
\quad \quad ,
\quad \quad
\deg (X_\cA) = \Vol(\cA)= \Vol_\cA (Q_\cA).
$$
Dans le cas o{\`u} 
$L_\cA = \Z^n$, la forme volume  $\Vol_\cA$ 
co{\"\i}ncide avec $n!$ fois la forme volume euclidienne 
$\Vol_{n}$ de $\R^n$; en particulier  
$\dim(X_{\cA}) = n$ et $\deg(X_{\cA}) = n!\, \Vol_{n}(Q_\cA)$.

\smallskip

Nous montrons que les minimums successifs de $X_\cA$ peuvent s'expliciter en termes de la 
combinatoire de l'ensemble $\cA$. 
Soit \,$\Faces(Q_\cA)$\,  
l'ensemble des faces du polytope $Q_\cA$, et 
pour $i=0, \dots, r:=\dim (\cA) $  posons 
$$
N_\cA(i) := \min \,\Big\{ \ \Card(P \cap \cA)  \ ; \
P \in \Faces(Q_\cA) , \ \dim (P) = i  \,\Big\}
$$
pour le minimum des cardinaux de l'ensemble  
$\cA$ restreint aux faces de $Q_\cA$ de dimension $i$.

\begin{thm} \label{1} 
\ $
\displaystyle{\mu_i (X_\cA) = \frac{1}{2}   \, \log N_\cA(r-i+1)}
$ \ pour \ $  i=1, \dots, r+1$.  
\end{thm}

En particulier  \,$N_\cA(r) = \Card(\cA) = N+1$\, et  \,$ N_\cA(0) =1$
et donc 
$$
\mu^\ess(X_\cA) = \frac{1}{2} \, \log (N+1)
\quad \quad , \quad \quad
\mu^\abs(X_\cA) = 0.
$$

Comme  corollaire du th{\'e}or{\`e}me~\ref{1} et 
du th{\'e}or{\`e}me des minimums successifs de Zhang,
on d{\'e}duit l'encadrement suivant pour la hauteur d'une vari{\'e}t{\'e} 
torique  (Corollaire~\ref{altura}):  
$$
\frac{1}{2} \, \log(N+1) \, \Vol (\cA)
\, \le \,
h(X_\cA)
\, \le \,
\frac{(r+1)}{2} \, \log(N+1) \, \Vol (\cA).
$$

{\`A} l'heure actuelle, on ne dispose  d'aucune expression exacte
(autre que sa d{\'e}finition)
pour $h(X_\cA)$, donc 
cet encadrement est  
significatif. 
En fait, il serait fort int{\'e}ressant de trouver 
une telle expression {exacte}  
pour la hauteur d'une vari{\'e}t{\'e} torique projective
dans le cas g{\'e}n{\'e}ral. 
Remarquons 
que pour les espaces projectifs  
(le seul cas des 
vari{\'e}t{\'e}s toriques 
dont on conna{\^\i}t la hauteur projective)
le th{\'e}or{\`e}me~\ref{1} montre qu'aucune des estimations dans le th{\'e}or{\`e}me
des minimums successifs de Zhang n'est 
exacte, {\it voir}~le paragraphe~\ref{Hauteur}.

A partir de ces estimations, on d{\'e}duit 
ais{\'e}ment la  majoration explicite suivante  
pour la hauteur du 
$\cA$-r{\'e}sultant $\Res_\cA$, 
dont on renvoie au paragraphe~\ref{Hauteur} pour la d{\'e}finition pr{\'e}cise. 
On note \,$ h_\mmax(\Res_\cA) $\, 
la hauteur associ{\'e}e {\`a} la norme sup, 
d{\'e}finie comme {\'e}tant le logarithme du maximum  des valeurs absolues 
des coefficients de $ \Res_\cA$.

\begin{cor} \label{res} 

Soit \,$\cA =\{ a_0, \dots, a_N \} \subset \Z^n$\, un ensemble fini 
tel que 
 \,$L_\cA = \Z^n$, alors 
$$
h_\mmax (\Res_\cA) 
\, \le \,    \frac{3}{2} \, (n+1) \, \log (N+1) \, \Vol(\cA).
$$
\end{cor} 

Les seules estimations
pour la hauteur de $\Res_\cA$
dont   
on disposait auparavant  
{\'e}taient celles qu'on peut obtenir {\`a} partir des algorithmes 
pour son  calcul~\cite[\S~7.6]{CoLiOS98},
\cite{DAndrea01}, {\it voir} aussi~\cite[Prop.~1.7]{KrPaSo01} ou encore~\cite[Thm. 23]{Rojas00}. 
Celles-ci  
sont  loin d'{\^e}tre pr{\'e}cises, et notre estimation les am{\'e}liore 
d'un facteur exponentiel.
Par ailleurs, la seule information 
dont on dispose
concernant  la valeur {\it exacte} 
des coefficients de $\Res_\cA$ est pour 
les coefficients extr{\'e}maux, 
{\'e}gaux {\`a} $\pm 1$~\cite[Cor.~3.1]{Sturmfels94}.

\medskip

La majoration de $h(X_\cA) $  
nous permet aussi d'estimer la 
hauteur des sous-vari{\'e}t{\'e}s d{\'e}finies par des {\'e}quations {\`a} support 
restreint.
Soit
\,$Q \subset \R^n$\, un polytope rationnel, c'est-{\`a}-dire 
dont les
sommets sont des vecteurs entiers. 
Pour un point 
$\xi \in \T^n$, on d{\'e}finit sa {\em $Q$-hauteur de Weil} 
par la formule
$$
\widehat{h}_Q(\xi) := 
\sum_{v \in M_K} 
\frac{[K_v : \Q_v]}{ [K:\Q]}  \, \log \max 
\Big\{ \,|\xi^a|_v \
; \ a \in Q \cap \Z^n
\,\Big\} 
\ \in \R_+,  
$$
o{\`u} $K$ est un corps de nombres  contenant les coordonn{\'e}es de
$\xi$ et $M_K$ est l'ensemble des places de $K$.

On pose
\,$\Vol_n(Q) $\,  
pour le volume euclidien de $Q$.
On pose aussi \ $
||f||_1 :=\sum_{a\in \Z^n} |c_a|$ pour la norme  $ \ell^1$ 
d'un polyn{\^o}me de Laurent 
\,$f  = \sum_{a\in \Z^n} c_a \, x^a 
\in \Z[x_1^{\pm 1} , \dots, x_n^{\pm 1} ]$.

\begin{thm} \label{2} 
Soient 
\,$f_1, \dots, f_n \in \Z[x_1^{\pm 1}, \dots, x_n^{\pm 1}] 
$\,
des polyn{\^o}mes de Laurent
{\`a} coefficients entiers, et posons 
\,$ Q:= \NP(f_1, \dots, f_n) \subset \R^n$\,  
pour leur polytope de Newton. 

Soit 
\,$Z(f_1,\dots, f_n)_0 $\, l'ensemble des points isol{\'e}s de 
\,$Z(f_1, \dots, f_n) \subset \T^n$. 
Pour chaque point   $\xi$ dans cet ensemble on note
\ $\ell(\xi) := \dim_\Qbarra  
\left(
\Qbarra[x_1^{\pm 1} , \dots, x_n^{\pm 1} ]
/ (f_1, \dots, f_n) \right)_{I(\xi)}$\ la multiplicit{\'e}
d'intersection de
$f_1, \dots, f_n$ en  $\xi$.   
Alors
$$
\sum_{\xi \in Z(f_1,\dots, f_n)_0}  \ell(\xi) \, \widehat{h}_Q (\xi) 
\ \le \ 
n!\, \Vol_n (Q)  \, \sum_{i=1}^n \, \log||f_i||_1. 
$$

\end{thm}

Ceci est un analogue arithm{\'e}tique du th{\'e}or{\`e}me classique
de B{\'e}zout-Koushnirenko~\cite[\S5.5]{Fulton93}, 
\cite[\S~6.2, Thm.~2.2]{GeKaZe94}, 
\cite[\S~7.5, Th.~5.4]{CoLiOS98}: 
$$
         \sum_{\xi \in Z(f_1,\dots, f_n)_0}   \ell(\xi)  
\ \le \ n! \, \Vol_n(Q) .  
$$ 
Ce r{\'e}sultat  est  un raffinement du th{\'e}or{\`e}me de B{\'e}zout
arithm{\'e}tique pour les hypersurfaces, 
et 
il am{\'e}liore  et rend compl{\`e}tement {e}ffectif
le cas   
non-mixte du th{\'e}or{\`e}me de
Bernstein-Koushnirenko arithm{\'e}tique d{\^u} {\`a}
V.~Maillot~\cite[Cor.~8.2.3]{Maillot00}, 
{\it voir} la remarque~\ref{rem-maillot}.

Sa d{\'e}monstration 
est bas{\'e}e sur la majoration de $h(X_\cA)$ d{\'e}j{\`a} mentionn{\'e}e, 
le th{\'e}or{\`e}me de B{\'e}zout
arithm{\'e}tique classique, et
des r{\'e}sultats plut{\^o}t {\'e}l{\'e}mentaires de la th{\'e}orie
de l'intersection g{\'e}om{\'e}trique. 
On obtient aussi une variante  pour les 
intersections arbitraires (Proposition~\ref{BK-gral}).  
Signalons qu'il serait d'un grand int{\'e}r{\^e}t d'{\'e}tendre ces r{\'e}sultats 
au cas mixte g{\'e}n{\'e}ral; 
nous esp{\'e}rons pouvoir le faire {\`a} l'aide d'une extension des techniques
introduites dans le pr{\'e}sent article.

\smallskip

Le $\cA$-r{\'e}sultant et le th{\'e}or{\`e}me de Koushnirenko 
sont {\`a} la base de la th{\'e}orie de 
l'{\'e}limination dite ``sparse'' (creuse), {\it voir}~\cite{CoLiOS98}, \cite{Sturmfels02}.  
Dans ce contexte,
il est peut-{\^e}tre int{\'e}ressant de signaler que 
 le th{\'e}or{\`e}me~\ref{2} 
montre  que 
la complexit{\'e} binaire d'une repr{\'e}sentation symbolique 
de 
\,$ Z(f_1,\dots, f_n)_0 $\,
est {\em polynomiale} en  
la complexit{\'e} binaire de $f_1, \dots, f_n$ 
et en le  volume (normalis{\'e} par $n!$) du polytope de Newton 
correspondant.

\medskip

Comme application de ces r{\'e}sultats,
nous d{\'e}duisons  une  am{\'e}lioration significative  
du Nullstellensatz arithm{\'e}tique creux d{\^u} {\`a}  
T. Krick, L.M. Pardo et l'auteur~\cite[Cor.~3]{KrPaSo01}:

\begin{thm} \label{nss} 
Soient \,$f_1,\ldots,f_s \in \Z[x_1,\ldots,x_n]$\, des 
polyn{\^o}mes 
{\`a} coefficients entiers  
sans z{\'e}ros
communs dans
$\C^n$.

Posons \ $d:= \max_i \deg (f_i) $,   
\ $h:= \max_i h_\mmax (f_i)$ et 
 \ $\cA : = \Supp \Big(1, x_1, \dots, x_n, f_1, \dots , f_s\Big) 
\subset \N^n$; 
alors il existe \,$ a\in \Z \setminus \{0\}$\, et \,$g_1,\ldots,g_s
\in \Z [x_1,\ldots,x_n]$\, tels que

\smallskip 

\begin{itemize}
\item[$\bullet$]  $a \ = \  {g}_1 \, {f}_1 + \cdots + {g}_s \, {f}_s $,

\medskip 

\item[$\bullet$]  $ \deg (g_i)  \ \le \   2\, n^2 \, d\, \Vol(\cA) $,

\medskip 

\item[$\bullet$]  $ h_\mmax(a) , h_\mmax (g_i) \ \leq \ 
2\, (n+1)^3\,  d\, \, \Vol(\cA) \, \Big(
h +\log s    +  14 \, (n+1) \, d\,
\log(d+1) \Big)$. 
\end{itemize}

\end{thm}

Ce r{\'e}sultat repr{\'e}sente une am{\'e}lioration 
d'un facteur exponentiel
par rapport aux estimations connues; 
la majoration pour les hauteurs 
devienne ainsi {\em polynomiale} 
en tous les param{\`e}tres concern{\'e}s.

\bigskip

% ----------- Agradecimientos -----------------------------------------------

\noindent {\bf Remerciements.} 
Je tiens {\`a} remercier chaleureusement 
Patrice Philippon 
pour des tr{\`e}s nombreuses discussions et {\'e}claircissements, notamment 
autour
du th{\'e}or{\`e}me de B{\'e}zout. 
J'ai eu des discussions int{\'e}ressantes  avec  Marc Chardin 
et Vincent Maillot. 
Je remercie Teresa Krick pour sa lecture attentive d'une 
version pr{\'e}liminaire de ce texte; 
ainsi qu'au rapporteur anonyme pour des nombreuses remarques
utiles.

Bernd Sturmfels m'a encourag{\'e} {\`a} r{\'e}fl{\'e}chir sur la hauteur du
$\cA$-r{\'e}sultant. 
Le pr{\'e}sent article  est en partie cons{\'e}quence de mes tentatives
pour r{\'e}soudre ce probl{\`e}me.

% ----------- Seccion 0,5 -----------------------------------------------

\typeout{Hauteur des points et des polynomes} 

\section{Hauteur des points et des polyn{\^o}mes}

\label{points}

\setcounter{equation}{0}
 
\renewcommand{\theequation}{\thesection.\arabic{equation}}

On note  par  $\Q$ le corps des nombres rationnels,
$K$ un corps de nombres et 
$\Qbarra$ la cl{\^o}ture alg{\'e}brique de $\Q$. 
On note  par 
$\T^n$ le tore alg{\'e}brique
et $\P^N$ l'espace projectif
sur $\Qbarra$, de dimension $n$ et $N$ respectivement.  

Pour chaque premier rationnel $p$  on note $ | \cdot |_p $ la
valeur
absolue $p$-adique sur $\Q$ telle que
$|p|_p=p^{-1}$; 
on note aussi  $|\cdot|_\infty$
ou simplement $ | \cdot | $
la valeur absolue standard. 
Celles-ci forment un ensemble complet de valeurs absolues sur
$\Q$:
on identifie l'ensemble $M_\Q$
de ces valeurs absolues {\`a} l'ensemble $\{ \infty, \, p
\, ; \, p \ \mbox{premier} \}  $.
Plus g{\'e}n{\'e}ralement, on d{\'e}signe par $M_K$ l'ensemble des valeurs absolues
de $K$ {\'e}tendant
les valeurs absolues de $M_\Q$, 
et
on note $M_K^\infty$ le sous-ensemble de $M_K$
des valeurs absolues
archim{\'e}diennes. 

On note  $\R$ le corps des nombres r{\'e}els et 
$\C$ le corps des nombres complexes; on pose
$\R_+$  l'ensemble des nombres r{\'e}els non-n{\'e}gatifs. 
On note par 
$\Z$ l'anneau des entiers rationnels,
et par  $\N$ et $\N^\times$ les entiers naturels avec et sans
0, respectivement. 
Pour $ b = (b_0, \dots, b_N) \in \Z^{N+1}$, on pose $\deg(b):=
b_0+\cdots +b_N \in \Z$ pour le degr{\'e} du mon{\^o}me $x^b$ associ{\'e}.

\bigskip 

Dans ce paragraphe nous rappelons les d{\'e}finitions et propri{\'e}t{\'e}s de
base
des diff{\'e}rents notions de hauteur pour les points (et plus g{\'e}n{\'e}ralement,
pour les 
sous-vari{\'e}t{\'e}s) et pour les polyn{\^o}mes. 

\medskip

Pour un point $\xi=(\xi_0:\cdots:\xi_N) \in \P^N$ on note
\,$ h(\xi) $\, son hauteur {\it projective}, d{\'e}finie par la
formule
$$
h(\xi) : =  \sum_{v \in M_K^\infty } {\lambda_v(K)}  \,
\log \sqrt{ |\xi_0|_v^2  + \cdots +
|\xi_N|_v^2 }   \hspace{2mm}  +
\hspace{-2mm}  \sum _{v \in M_K \setminus M_K^\infty}
\lambda_v(K) \,  \log \max \Big\{ |\xi_0|_v , \dots,
|\xi_N|_v\Big\}, 
$$
o{\`u}   \,$K$\, est  un
corps de nombres contenant les $\xi_i$, et dont pour chaque 
pour $v \in M_K$ on pose 
$$
\lambda_v(K) := \frac{[K_v : \Q_v]}{[K:\Q]}. 
$$
Cette expression ne  d{\'e}pend pas  
du choix des coordonn{\'e}es
homog{\`e}nes de $\xi$, gr{\^a}ce {\`a} la formule du produit.

\smallskip 

Cette notion s'{\'e}tend au vari{\'e}t{\'e}s projectives de dimension sup{\'e}rieure,
suivant une construction due {\`a} P.~Philippon. 
Soit
$f = \sum_a c_a \, U^a \in K[U_0, \dots,  U_n]$ un polyn{\^o}me
en $n+1$ 
groupes $U_i$ de $N_i+1$ variables chacun.
Posons 
$S_{N_i+1}:= \Big\{ (z_0, \dots, z_{N_i}) \in \C^{N_i+1}
\, ; \ |z_0|^2 + \cdots + |z_{N_i}|^2 = 1 \Big\} $
la sph{\`e}re unit{\'e} de $\C^{N_i+1}$,
{\'e}quip{\'e}e de
la mesure
$\mu_{N_i+1}$
de masse totale 1 et invariante par rapport au groupe
unitaire $U(N_i+1)$.
Pour une place archim{\'e}dienne $v \in M_K^\infty$,
la {\it $S_{N_0+1} \times \cdots \times S_{N_n+1}$-mesure} 
de $f$ relative  {\`a} $v$
est d{\'e}finie par l'int{\'e}grale 
$$
m_v( f;S_{N_0+1} \times \cdots \times S_{N_n+1})
:= \int_{S_{N_0+1} \times \cdots \times S_{N_n+1}}
\log |f|_v
\ \mu_{N_0+1} \times \cdots \times \mu_{N_n+1} \enspace.
$$
Notons que dans le cas o{\`u} chaque groupe se r{\'e}duit {\`a} une seule variable 
$U_i:= \{u_i\}$, cette notion se sp{\'e}cialise en la {\it mesure de Mahler}
de $f$, c'est-{\`a}-dire 
$$
m(f; S_1^{n+1})= 
m(f):= \int_0^1 \cdots \int_0^1 \log\Big|f(e^{2\, \pi \, i \, u_0}, 
\dots, e^{2\, \pi \, i \, u_n}) \Big| \, du_0 \cdots du_n.
$$
Pour une place ultram{\'e}trique $ v \in M_K \setminus M_K^\infty$, on 
note $
|f|_v:= \max \Big\{ |c_a|_v \, ; \ a \in \N^{N_0+1}
\times \cdots \times \N^{N_n+1} \Big\}$ 
la {\it norme sup} de $f$ relative {\`a} $v$.

\smallskip

Maintenant soit $V \subset \P^N$
une $K$-vari{\'e}t{\'e} quasi-projective {\'e}quidimensionnelle de dimension $n$. 
Soit $\Ch_V \in K[U_0, \dots, U_n]$ sa {forme de Chow}, qui est
un polyn{\^o}me
en $n+1$ groupes de $N+1$ variables chacun,
homog{\`e}ne de degr{\'e} $\deg(V) $ en chaque groupe $U_i$.  
Suivant~\cite{Philippon919495}, la {\it hauteur projective} de $V$ est 
par d{\'e}finition 
\begin{eqnarray*}
h(V) &:=&
\sum_{v \in M_K^\infty} \lambda_v(K) \,
m_v(\Ch_V ; S_{N+1}^{n+1}) +
\sum_{v \in M_K \setminus M_K^\infty} \lambda_v(K) \,
\log |\Ch_V|_v 
\\[-2mm]
&&  + (n+1)\,
\bigg(\sum_{j=1}^{N} \frac{1}{2\, j} \bigg) \, \deg(V) \ \ \in \R_+.
\end{eqnarray*}
Alternativement, cette hauteur peut se d{\'e}finir {\it via} la
th{\'e}orie d'Arakelov
comme la hauteur $h_{\ov{\cO(1)}}(\Sigma)$
relative au fibr{\'e} en droites universel
$\ov{\cO(1)}$ muni de la m{\'e}trique
de Fubini-Study, 
de la cl{\^o}ture de Zariski $\Sigma$ de $V$
dans $\P^N_{\cO_K}$, {\it voir}~\cite[\S~3.1.3]{BoGiSo94}.

\smallskip

On consid{\'e}rera aussi la {\it hauteur de Weil} d'un point $\xi \in \P^N$,
d{\'e}finie par 
$$
\hnorm(\xi) : =  \sum_{v \in M_K}\lambda_v(K) \,  \log \max \Big\{ |\xi_0|_v , \dots,
|\xi_N|_v\Big\}. 
$$
Ceci se compare avec la hauteur projective; on a 
$$
\hnorm(\xi) \le h(\xi) \le \hnorm (\xi) + \frac{1}{2} \, \log(N+1).
$$
Le tore $\T^N$ peut s'identifier {\`a} l'ouvert $(\P^N)^\circ := \P^N \setminus
\{ (x_0:\cdots:x_N)\,;\ x_0 \cdots x_N =0\}$ {\it via} l'inclusion
$(t_1, \dots,t_N) \mapsto (1:t_1:\cdots:t_N)$. 
Ceci induit une notion de hauteur projective $h_\iota$ et de hauteur
de Weil $\hnorm_\iota$
pour les sous-vari{\'e}t{\'e}s et pour les points de $\T^N$, respectivement.

Plus g{\'e}n{\'e}ralement, 
{\'e}tant donn{\'e}  un ensemble fini
$\cA=\{ a_0, \dots, a_N\} \subset \Z^n$ 
on peut d{\'e}finir une notion de {\it $\cA$-hauteur
projective} pour une sous-vari{\'e}t{\'e} {\'e}quidimensionnelle $V \subset \T^n$ en posant 
$$
h_\cA(V) := h(\varphi_\cA(V)).
$$
Ceci se sp{\'e}cialise en la hauteur projective $h_\iota$ pour $\cA:= \{ 0, e_1,\dots,
e_n\}$, o{\`u} $e_j$ d{\'e}signe le
$j$-{\`e}me vecteur de la base standard de $\R^n$. 
Similairement on d{\'e}finit
la {\it $\cA$-hauteur de Weil} d'un  point $\xi \in \T^n$ par
$$
\widehat{h}_\cA(\xi) :=  \hnorm(\varphi_\cA(\xi));
$$
ceci se sp{\'e}cialise en la hauteur de Weil  $\hnorm_\iota$ pour $\cA:= \{ 0, e_1,\dots,
e_n\}$.

Remarquons que pour un polytope rationnel  \,$Q \subset \R^n$, 
la  {\it $Q$-hauteur de Weil} $\hnorm_Q$ 
co{\"\i}ncide avec la $\cA_Q$-hauteur de Weil,
pour  
\,$\cA_Q:= Q \cap \Z^n  \subset \Z^n$\, 
l'ensemble des points entiers dans $Q$.

\medskip 

On fera appel {\`a} plusieurs notions de hauteur pour les polyn{\^o}mes, 
suivant le choix d'une m{\'e}trique pour les places {\`a} l'infini. 
Soit $d \in \N$ et  $||\cdot||$ une m{\'e}trique sur la partie gradu{\'e}e 
$\C[x_0, \dots, x_N]_d$. 
Soit $v \in M_K^\infty$ une place archim{\'e}dienne 
et $\sigma_v:K \hookrightarrow \C$ une immersion correspondant {\`a} $v$. 
Pour un polyn{\^o}me homog{\`e}ne 
$f\in K[x_0, \dots, x_N]_d$ on
pose   alors 
$$
||f||_v:= ||\sigma_v(f)||
$$
pour la {\it norme} de $f$ relative {\`a} $v$. 
Rappelons que pour $v \in M_K \setminus M_K^\infty$ on note $|f|_v$
la norme sup de $f$ relative {\`a} $v$. 
Ainsi, la m{\'e}trique $||\cdot||$ d{\'e}finit une hauteur par la formule
$$
h_{||\cdot||}(f) := \sum_{v\in M^\infty}\lambda_v(K) \, \log||f||_v 
\ + \sum_{M_K\setminus M_K^\infty} \lambda_v(K) \, \log|f|_v. 
$$
Pour la norme sup et les normes $\ell^1$ et $\ell^2$,  on notera les
hauteurs correspondantes par $h\mmax$, $h_1$ et $h_2$, respectivement. 
Les relations entre ces normes fournissent des encadrements pour les
hauteurs correspondantes; on a 
$$
h_\mmax(\xi) \le h_2(\xi) \le h_1(\xi)
\quad \quad , \quad \quad 
h_1(\xi) \le h_2(\xi) + \frac{1}{2} \, \log {N+d \choose N} \le 
 h_\mmax(\xi) + \log{N+d \choose n}.
$$
Les hauteurs  $h\mmax$, $h_1$ et $h_2$ s'{\'e}tendent de fa{\c c}on naturelle 
au cas des polyn{\^o}mes pas forc{\'e}ment homog{\`e}nes, voire des polyn{\^o}mes de Laurent. 

\smallskip 

Aussi on consid{\'e}rera
la m{\'e}trique  \,$||\cdot ||_\W$\,  sur $\C[x_0, \dots, x_N]_d$ 
d{\'e}finie par
$$
|| f ||_\W^2 := 
\sum_{a; \, \deg(a)=d}  
{  {d \choose a}^{-1} \, |c_a|^2}
$$
pour
\,$f = \sum_a c_a \, x^a 
\in \C [x_0, \dots, x_N]_d$.
Pour \,$\xi \in \C^{N+1}$, l'in{\'e}galit{\'e} de Cauchy-Schwartz entra{\^\i}ne 
\begin{equation} \label{weyl} 
|f(\xi)|  = \left| \sum_a c_a \, \xi^a \right| 
\le \bigg( \sum_{a}  
{  {d \choose a}^{-1} \, |c_a|^2}
\bigg)^{1/2}  \, \bigg(\sum_a {d \choose a }
  |\xi^a|^2 \bigg)^{1/2} 
 = ||f||_\W \, ||\xi||_2^d.
\end{equation} 

Pour  un polyn{\^o}me homog{\`e}ne  $f \in K[x_0, \dots, x_N]_d$, 
on notera respectivement $||\cdot||_{\W,v}$ et $h_\W$ la m{\'e}trique 
relative {\`a} une place $v \in M_K$ et 
la hauteur correspondantes {\`a} cette m{\'e}trique. 
L'indice  est pour 
H.~Weyl, le premier {\`a} notre  connaissance   
{\`a} consid{\'e}rer cette m{\'e}trique, {\it voir}~\cite[\S~3.7]{Weyl39}.

% ----------- Seccion 1 -----------------------------------------------

\typeout{Minimums successifs} 

\section{Minimums successifs}

\label{minimums}

Soit
\,$\Tors(\T^N)$\, le sous-groupe des points de torsion de $\T^N$,
on a $\xi=(\xi_1,\dots,x_N) \in \Tors(\T^N)$ si et seulement 
si $\xi_i$ est une racine de l'unit{\'e} pour $i=1,\dots, N$.
Le lemme suivant montre que 
le minimum de la hauteur projective $h_\iota$ sur les points de $\T^N$ est atteint aux 
points de torsion.

\begin{lem} \label{torsion}
Soit  $\xi \in \T^N $; alors
\ $
\displaystyle{
h_\iota(\xi) \ge \frac{1}{2}  \, \log (N+1)} 
$,  
avec {\'e}galit{\'e}  si et seulement si \  $\xi \in \Tors(\T^N)$.

\end{lem}

\begin{demo}
Soit $ v \in M_K^\infty$.
L'in{\'e}galit{\'e} arithm{\'e}tico-g{\'e}om{\'e}trique~\cite[Thm.~9]{HaLiPo88}
entra{\^\i}ne
$$
\frac{1+|\xi_1|_v^2+ \cdots + |\xi_N|_v^2}{N+1}
\ge  \Big(|\xi_1|_v \cdots
|\xi_N|_v \Big)^{2/(N+1)}
$$
avec {\'e}galit{\'e} si et seulement si \ $
|\xi_i |_v = 1 $
pour tout $i$. 
{\'E}quivalemment
\begin{equation} \label{in1}
 \log  \sqrt{1+|\xi_1|_v^2+ \cdots +
|\xi_N|_v^2} 
 \ge
 \frac{1}{2} \, \log (N+1) +
\frac{1}{N+1} \, \log \Big(|\xi_1|_v \cdots |\xi_N|_v\Big).
\end{equation}
Pour $ v \in M_K \setminus M_K^\infty$ on a
\begin{equation} \label{in2}
\log \max \Big\{ 1, |\xi_1|_v , \dots,
|\xi_N|_v \Big\}  \ge
\frac{1}{N+1}
 \, \log \Big( |\xi_1|_v \cdots |\xi_N|_v \Big) ,
\end{equation}
avec {\'e}galit{\'e} si et seulement si \ $
|\xi_i |_v = 1 $
pour tout $i$, donc
\begin{eqnarray*}
h_\iota(\xi) & = & \sum_{v \in M_K^\infty} {\lambda_v(K)}  \,
\log  \sqrt{1+|\xi_1|_v^2+ \cdots +
|\xi_N|_v^2} \hspace{3mm}   
+ \hspace{-4mm} \sum_{v \in M_K \setminus M_K^\infty}\hspace{-3mm}
\lambda_v(K) \,  \log \max \Big\{ 1, |\xi_1|_v , \dots,
|\xi_N|_v \Big\} \\[2mm]
& \ge &
\frac{1}{2} \, \log (N+1) +
\frac{1}{N+1} \,
\sum_{v \in M_K}
\lambda_v(K)
 \, \log \Big(|\xi_1|_v \cdots |\xi_N|_v \Big) =
\frac{1}{2} \, \log (N+1)
\end{eqnarray*}
par la formule du produit.
Pour avoir l'{\'e}galit{\'e}, il faut que 
les in{\'e}galit{\'e}s~(\ref{in1}) et~(\ref{in2}) soient des {\'e}galit{\'e}s. 
Ceci {\'e}quivaut {\`a} ce que
\,$
|\xi_i |_v = 1 $\,
pour tout $v \in M_K$ 
et $i=1, \dots, N$,
c'est-{\`a}-dire  $ \xi \in \Tors(\T^N)$.
\end{demo}

Consid{\'e}rons l'action diagonale de $\T^n$ sur $\P^N$ associ{\'e}e {\`a}
l'ensemble $\cA$ 
$$
*_\cA: \T^n \times \P^N \to \P^N
\quad \quad , \quad \quad
(t, x) \mapsto t *_\cA x:= (t^{a_0}\, x_0 : \cdots : t^{a_N} \, x_N).
$$
Les orbites de l'action restreinte {\`a} $X_\cA$ 
sont en correspondance avec l'ensemble $\Faces(Q_\cA)$
des faces du polytope $Q_\cA$: 
pour chaque face $ P$ on consid{\`e}re un point 
$e_P := (e_{P, \, 0} : \cdots : e_{P, \, N} ) \in \P^N$ 
d{\'e}fini par $  e_{P, \, j}  :=1$ 
si \,$ a_j \in P $  et  $  e_{P, \, j}
:= 0$ sinon; 
la bijection est  donn{\'e}e par~\cite[Ch.~5, Prop.~1.9]{GeKaZe94},  
\cite[\S~3.1]{Fulton93} 
$$
P \mapsto X_{\cA, P}^\circ := \T^n *_\cA e_P \ \subset \P^N 
\enspace. 
$$
Ainsi la d{\'e}composition de $X_\cA$ comme union
disjointe des orbites de cette action s'explicite en
\begin{equation} \label{orbitas} 
X_\cA  = \bigsqcup_{P \in F(Q_\cA) } X_{\cA, P}^\circ.
\end{equation}
Posons  $N(P):= {\rm Card}\{ i \, ; \ a_i \in P\} -1$ et
$
\cA(P) = (a_i \, ; \ a_i \in P) \in(\Z^n)^{N(P)+1}$.
On v{\'e}rifie que $X_{\cA, P}^\circ \subset \P^N$ 
est l'orbite principale d'une vari{\'e}t{\'e} torique 
contenue
dans un sous-espace standard $E \cong \P^{N(P) }$. 
Cette vari{\'e}t{\'e} torique 
s'identifie {\`a} la sous-vari{\'e}t{\'e} 
$ X_{\cA(P)}^\circ \subset \P^{N(P) }$ 
de dimension {\'e}gale {\`a} la dimension (r{\'e}elle) de la face  $P$, 
{\it via} l'inclusion 
canonique $i_P: \P^{N(P) } \hookrightarrow \P^N$.  
On renvoie le lecteur int{\'e}ress{\'e} {\`a}~\cite{Fulton93}, \cite{GeKaZe94} ou
\cite{CoLiOS98} 
\ pour les propri{\'e}t{\'e}s de base des vari{\'e}t{\'e}s toriques projectives.

Le th{\'e}or{\`e}me~\ref{1} est   cons{\'e}quence directe 
de la d{\'e}composition ci-dessus 
et du
lemme~\ref{torsion}. 
Montrons d'abord quelques 
propri{\'e}t{\'e}s {\'e}l{\'e}mentaires des minimums successifs.

\begin{lem} \label{union}
\ 

\begin{itemize} 

\medskip 

\item[(a)] 
Soit \,$V \subset \P^N $\, 
une vari{\'e}t{\'e} quasi-projective quelconque, 
et \,$ Z \subset V$\, une sous-vari{\'e}t{\'e} quasi-projective 
de dimension $s$, alors
\ $
\mu_i(V) \le \mu_{s-r+i} (Z) 
$ \ 
pour $ i =r-s+1, \dots , r+1$. 

\medskip 

\item[(b)]
Soit \,$ V =   \bigcup_{j \in J } Z_j $\, un recouvrement fini 
de $V$ par des sous-vari{\'e}t{\'e}s 
quasi-projectives de dimension \  $s_j := \dim (Z_j) $, 
alors pour $ i =1, \dots, r+1$
$$
\mu_i (V) = \min \, \Big\{ \, \mu_{ s_j- r +i } (Z_j) \ ; 
\ j \in J, \ s_j \ge r-i+1 \Big\}. 
$$
\end{itemize} 
\end{lem}

\begin{demo} 
Partie (a):  
Posons $\ell:= r-i$. 
Pour  $\varepsilon > 0$ donn{\'e}, on prend une 
sous-vari{\'e}t{\'e}  
$W \subset V$ de codimension $i$ (ou {\'e}quivalemment
telle que \,$\dim ( W) = \ell$)
suffisamment grande pour que 
\ $
\mu_i (V) \le \mu^\abs (V \setminus W) + \varepsilon$. 

On a   \,$\dim ( W \cap Z ) \le \ell$\, et donc \ 
$
\mu^\abs (Z \setminus W ) \le \mu_{s-r+i}  (Z) 
$ \ 
par d{\'e}finition, car 
\,$ \codim_{Z} ( W \cap Z ) \ge s -\ell = s-r+i$. 
Puis 
$$ \mu^\abs (V \setminus W) 
\le \mu^\abs (Z \setminus W )
$$ 
car \,$V \setminus W  \supset Z \setminus W$; 
on en conclut \ $ \mu_i (V) \le \mu_{s-r+i}  (Z) $.

\medskip

Partie (b): Par application directe de la partie (a), 
on voit 
que 
le premier terme est born{\'e} par le deuxi{\`e}me. 
Ainsi on se r{\'e}duit {\`a} d{\'e}montrer l'autre in{\'e}galit{\'e}. 

\smallskip 

Posons $\ell:= r-i$.  
Pour chaque $ j \in J$ on prend une sous-vari{\'e}t{\'e} \,$W_j \subset
Z_j$\, 
de la fa{\c c}on suivante: dans le cas \,$s_j \ge \ell+1$\, on prend $W_j$  de dimension
$\ell$, suffisamment grande pour que  
\ $
\mu_{s_j-r+i} (Z_j) \le \mu^\abs (Z_j \setminus W_j ) + \varepsilon, 
$  
tandis que pour le cas \,$s_j \le \ell$\, on prend 
\ $W_j:= Z_j$. 

Soit 
$$ 
W :=  \overline{\bigcup_{j\in J} W_j} \, \cap Z  \  \subset Z 
$$  
la cl{\^o}ture de Zariski dans $Z$    de la r{\'e}union des $W_j$. 
Donc \,$\dim ( W) \le \ell$, 
 car $J$ est un 
ensemble
 fini et    $\dim (W_j) \le \ell$ pour tout $j \in J$. 
Ainsi 
\ $
\mu_i(V) \ge \mu^\abs(V \setminus W)$ et de plus 
$$
 \mu^\abs(V \setminus W) \, \ge \,  \min \, \{ \, \mu^\abs 
 (Z_j \setminus W_j ) \ ; 
\ j \in J, \ s_j \ge r-i+1 \}
$$
car 
 \ $
V \setminus W \subset   
\bigcup_{ \{ j \, ; \, s_j \ge \ell+1 \} } (Z_j \setminus W_j)$. 
On en d{\'e}duit 
\ $
\mu_i(V) \ge 
 \min \, \{ \, \mu_{s_j-r+i} (Z_j) \, ; 
\ j \in J, \ s_j \ge r-i+1 \} - \varepsilon 
$, 
ce qui {\'e}tablit 
l'{\'e}galit{\'e}.
\end{demo}

\begin{Demo}{D{\'e}monstration du th{\'e}or{\`e}me~\ref{1}.--} 
Soit \,$\cB \subset \Z^n$\, un ensemble fini quelconque, 
et posons  
\,$M+1:= \Card(\cB) $\, pour son cardinal et  
 \ $
X_\cB^\circ := \varphi_\cB(\T^n) \subset \P^M 
$\ 
pour le tore associ{\'e}. 
Ce tore  
$ X_\cB^\circ$ est en fait 
contenu dans \ $(\P^M)^\circ$ et par le lemme~\ref{torsion} on a 
$$ 
\mu^\abs ( X_\cB^\circ) \ge \frac{1}{2} \, \log ( M+1). 
$$

En outre \,$\Tors(\T^n)$\, est un ensemble  
dense de $\T^n$, et donc 
 $ \varphi_\cB(\Tors(\T^n))$ est aussi  un ensemble  dense de  
$ X_\cB^\circ$. 
On a \,$ h (\xi ) =  \frac{1}{2} \, \log ( M+1)$\, pour tout 
$ \xi$ dans  cet ensemble, car il est contenu dans $\Tors(\T^M)$. 
Alors
$$ 
\mu^\ess ( X_\cB^\circ) \le \frac{1}{2} \, \log ( M+1) , 
$$ 
et donc \ $
\mu_j(X_\cB^\circ)  = 
\frac{1}{2}  \, \log ( M+1)  
$ \ 
pour tout $1 \le j \le \dim (X_\cB) +1 $, {\`a} cause de l'encadrement
\,$ \mu^\abs(X_\cB^\circ) 
\le \mu_j(X_\cB^\circ) \le
\mu^\ess(X_\cB^\circ)$. 

\smallskip 
 
Maintenant soit \,$P \in \Faces(Q_\cA)$,  et soit 
 \ $X^\circ_{\cA,P} = \T^n *_\cA e_P \subset \P^N$ l'orbite torique
 associ{\'e}e.   
La hauteur des points est invariante par l'inclusion $i_P$ 
et donc 
$$
\mu_i \Big(X^\circ_{\cA,P}\Big)  = \mu_i \Big(  X_{\cA(P)}^\circ \Big) = 
 \frac{1}{2} \, \log N(P) 
$$
pour $i =1, \dots, \dim (P) +1$. 
Par le lemme~\ref{union}(b) 
et la d{\'e}composition~(\ref{orbitas}) 
on conclut 
$$
\mu_i (X_\cA) = 
\min \Big\{ \, \mu_{\dim(P) - r +i} \Big(X^\circ_{\cA,P}\Big) \ ; \ 
\dim (P) \ge r-i+1 \,\Big\}  
=
\frac{1}{2} \, \log N_\cA(r-i+1), 
$$
car ce minimum est atteint sur l'ensemble des faces de $Q_\cA$ 
de dimension $r-i+1$.
\end{Demo}

Comme on l'a d{\'e}j{\`a} remarqu{\'e} dans l'introduction, en particulier on a 
\,$\displaystyle \mu^\ess(X_\cA) = \frac{1}{2}   \, \log (N+1) $\, et \,$ 
\mu^\abs(X_\cA) = 0$. 

\smallskip 

Ce r{\'e}sultat montre  
aussi que l'ensemble des minimums successifs de $X_\cA$ peut
{\^e}tre tr{\`e}s vari{\'e}, 
suivant la combinatoire de l'ensemble $\cA$ (sauf pour  les conditions 
\ $ \mu_{1} (X_\cA) \ge \cdots \ge \mu_{r+1} (X_\cA) =0 $). 
Voici quelques exemples:

\medskip 

\begin{itemize}

\item[$\bullet$] L'espace projectif  \,$\P^n$\,  correspond {\`a}
  l'ensemble 
\,$\{ 0, e_1, \dots,  e_n \}  \subset \Z^n$ ($e_j$ est le
$j$-{\`e}me vecteur de la base standard de $\R^n$) donc
$$
\mu_i(\P^n) = \frac{1}{2} \, \log (n-i+2)
\quad \quad , \quad \quad 
i =1, \dots, n+1.
$$

\item[$\bullet$] Soit
\,$
V_{n,\, d} \subset \P^{{d+n \choose n} -1}
$ 
la vari{\'e}t{\'e} de Veronese, d{\'e}finie comme 
l'image de l'application $\P^n \to \P^{{d+n \choose n} -1}$, 
$x \mapsto \Big( x^b \, ; \ b \in \N^{n+1}, 
\ \deg(b) = d \Big)$. 
Ceci correspond {\`a} l'ensemble
\,$ \{ a \in \N^n \, ; \, \deg(a) \le d \} \subset \Z^n $\, et
donc
$$
\mu_i(V_{n, \, d}  ) = \frac{1}{2} \, \log { d+ n -i+1 \choose n-i+1} 
\quad \quad , \quad \quad 
i =1, \dots, n+1.
$$

\item[$\bullet$] Soit  \,$S_n \subset \P^{2^n-1} $\, la vari{\'e}t{\'e} de
  Segre, d{\'e}finie comme l'image de la immersion  
 $
(\P^1)^n
  \hookrightarrow \P^{2^n-1} $ {\it via} 
$((x_{10}:x_{11}), \dots, (x_{n0}:x_{n1}))
\mapsto \Big( \prod_{i=1}^n x_{i\, j_i} \, ; \ j_i \in \{0,1\}\Big)$. 
Ceci correspond {\`a}
  l'ensemble
\,$ \{ 0, 1\}^n \subset \Z^n$\, et donc 
$$
\mu_i (S_n) = \frac{\log 2}{2} \, (n-i+1) 
\quad \quad , \quad \quad 
i =1, \dots, n+1.
$$
\end{itemize}

\bigskip

Une {\em sous-vari{\'e}t{\'e} de torsion} de $\T^N$ est par d{\'e}finition
une sous-vari{\'e}t{\'e} de
la forme \,$\omega \cdot H$, o{\`u} $H$ est un sous-groupe alg{\'e}brique de
$\T^N$ et $\omega \in \Tors(\T^N)$~\cite[\S~6]{Zhang95}, 
\cite[\S~5]{Bilu97}. 
La proposition~\ref{bogomolov} ci-dessous 
montre que les
sous-vari{\'e}t{\'e}s de torsion r{\'e}alisent, pour la hauteur projective, 
les minimums
successifs les plus petits possibles
parmi les sous-vari{\'e}t{\'e}s de $\T^N$.

\begin{lem}

Soit \,$H \subset \T^N$\,  un groupe alg{\'e}brique de dimension $r$, 
alors il existe un ensemble
fini \,$ \cA = \{ a_0, \dots , a_N\} \subset \Z^r $\, 
et un sous-groupe 
fini \,$G \subset \T^N$\,  
tels que \ 
$
\displaystyle 
H = G \cdot X_\cA^\circ = \cup_{\omega \in G} \omega \cdot X_\cA^\circ$. 
\end{lem}

\begin{demo}
Notons d'abord que gr{\^a}ce {\`a}~\cite[\S~3.2.3, Thm.~5]{OnVi90}
il existe un isomorphisme
 \,$ \psi : \T^N \to \T^N$\, et des entiers positifs 
\,$c_{r+1}, \dots, c_N \in \N^\times$\, tels que
$$
\psi (H) = Z(y_{r+1}^{c_{r+1}}-1 , \dots, y_N^{c_N} -1)   \ 
\subset \T^N. 
$$
Alors on a \,$ \psi (H) = F \cdot K$ avec  
\,$ K:=  Z( y_{r+1}-1 , \dots, y_N -1)  \subset \T^N$\, et  
$$
F:= Z( y_1-1, \dots, \, y_r-1,  
\, y_{r+1}^{c_{r+1}} - 1 , \dots, \, y_N^{c_N} - 1)  \ \subset \T^N. 
$$

Par 
\cite[\S~3.2.3, Thm.~4 et Prob.~10]{OnVi90}
il existe des vecteurs entiers  $v_1, \dots , v_N \in \Z^N$ tels que 
l'inverse \,$\psi^{-1} : \T^N \to \T^N$\, s'{\'e}crit
comme \,$ y \mapsto (y^{v_1}, \dots, y^{v_N})$.  

\smallskip 

Posons 
 \,$a_i := ( v_{i \, 1} ,\dots, v_{i\, r} ) \in \Z^r$\, pour $i =1,
\dots, N$; ainsi on a   \,$\psi^{-1} (K) = \varphi_\cA(\T^r) =
 X_\cA^\circ $, \,$G:= \psi^{-1} (F)$\, 
est un groupe fini et
$$
H= \psi^{-1} (F \cdot K)  
=\psi^{-1} (F) \cdot \psi^{-1} (K) 
= G \cdot X_\cA^\circ. 
$$ 
\end{demo} 

\vspace{-5mm} 

\begin{prop} \label{bogomolov}  
\ 

\begin{itemize} 

\smallskip 

\item[(a)] 
Soit $V \subset \T^N $, alors
\ $ 
\displaystyle{ \mu_i (V) \ge \frac{1}{2}   \, \log (N+1)} 
$ \ 
pour  \ $ 
i =1, \dots, \dim (V) +1$. 

%\smallskip 

\item[(b)] 
Soit $V \subset \T^N  $ une sous-vari{\'e}t{\'e} de torsion, alors
\ $
\displaystyle{\mu_i (V) = \frac{1}{2}  \, \log (N+1) } $ \ pour
\ $i =1, \dots, \dim (V) +1$.

\end{itemize} 

\end{prop} 

\begin{demo}
La partie (a) est cons{\'e}quence directe  du  lemme~\ref{torsion}. 
Pour la partie (b), on remarque que le lemme pr{\'e}c{\'e}dent implique qu'il 
existe des  ensembles
finis \, $ \cA = \{ a_0, \dots , a_N\} $ $\subset \Z^r $ \, 
et 
\, $J \subset  \Tors(\T^n)$ \, tels que  
\,$
V = J \cdot X^\circ_\cA
$. 
Par le lemme~\ref{union}(b), il suffit  de consid{\'e}rer le cas 
irr{\'e}ductible, c'est-{\`a}-dire \, $V = \omega \cdot X^\circ_\cA$ \, 
avec \,  $\omega \in \Tors(\T^N) $.

\smallskip

L'ensemble \,$\omega \cdot \varphi_\cA (\Tors(\T^r)) \subset \T^N  $\, 
est dense dans \,$ \omega \cdot
X_\cA^\circ$\, et 
il est contenu dans $\Tors(\T^N)$, et donc
  $ h_\iota(\xi  ) = \frac{1}{2} \, \log (N+1)$ pour tout point $\xi
$ de cet ensemble.
On en conclut que  
\,$  \mu_i(V) \le \mu^\ess (V) \le \frac{1}{2} \, \log(N+1) $\, 
pour tout $i$, ce qui ach{\`e}ve la
d{\'e}monstration. 
\end{demo}

En particulier \,$
\displaystyle \mu_i(\T^N) = \frac{1}{2}  \, \log(N+1) $\, pour 
$i=1, \dots, N+1$. 

\smallskip 

En vue de ce r{\'e}sultat, il est naturel de se demander si la
propri{\'e}t{\'e} (b) caract{\'e}rise les sous-vari{\'e}t{\'e}s de torsion, parmi toutes les
sous-vari{\'e}t{\'e}s irr{\'e}ductibles de $\T^N$. 
Plus pr{\'e}cis{\'e}ment: soit $V \subset \T^N$ une sous-vari{\'e}t{\'e} 
irr{\'e}ductible qui n'est pas de torsion,
et posons  \,$V^\circ$\, 
pour la r{\'e}union des sous-vari{\'e}t{\'e}s de torsion
contenues dans $V$ 
et \,$V^* := V \setminus V^\circ$.
Est-il vrai que  
\begin{equation}  \label{X} 
\mu^\abs ( V^*) \, > \,   \frac{1}{2} \, \log (N+1) \ ? 
\end{equation} 

Par le th{\'e}or{\`e}me de M.~Laurent~\cite[Thm.~2]{Laurent84}, il n'existe qu'un 
nombre fini de  sous-vari{\'e}t{\'e}s de torsion maximales
contenues dans $V$, ce qui implique 
que 
$V^*$
est un ouvert non vide.
Ainsi 
\,$
\mu^\ess(V) \ge \mu^\abs(V^*) $ et 
donc une r{\'e}ponse affirmative {\`a} la question~(\ref{X}) impliquerait une
r{\'e}ponse aussi affirmative {\`a} la question soulev{\'e}e pr{\'e}c{\'e}demment, 
c'est-{\`a}-dire au fait que les sous-vari{\'e}t{\'e}s irr{\'e}ductibles de $\T^N$
soient caract{\'e}ris{\'e}es  par la propri{\'e}t{\'e} (b) ci-dessus.

\smallskip

Prenons {\`a} pr{\'e}sent la notation dans~\cite[\S~6]{Zhang95}. 
Soit \,$e_{V^*}$\, 
le minimum absolu de $V^*$ par rapport
{\`a} 
la hauteur de Weil \,$\widehat{h}_\iota$. 
Pour un point quelconque 
$\xi \in \T^N$ on a \ 
$\displaystyle \widehat{h}_\iota (\xi ) + \frac{1}{2}  \, \log (N+1) \ge h_\iota(\xi) $ \ et donc 
$$
e_{V^*} 
+ \frac{1}{2} \, \log (N+1) \ge 
\mu^\abs(V^*).  
$$
On en d{\'e}duit que une r{\'e}ponse affirmative {\`a} la question~(\ref{X}) ci-dessus
repr{\'e}senterait aussi  un raffinement du probl{\`e}me  de
Bogomolov sur le tore $\T^N$~\cite[Thm.~6.2]{Zhang95}, \cite[Thm.~5.1(b)]{Bilu97}.

% ----------- Seccion 2 -----------------------------------------------

\typeout{Hauteur}

\section{Estimations des hauteurs}

\label{Hauteur}

\setcounter{equation}{0}
 
Dans ce paragraphe nous obtenons 
les estimations pour la
hauteur  des vari{\'e}t{\'e}s projectives 
{\`a} partir des r{\'e}sultats du paragraphe pr{\'e}c{\'e}dent.

\bigskip 

La r{\'e}partition  des points alg{\'e}briques de petite hauteur d'une
sous-vari{\'e}t{\'e} est en {\'e}troite relation  
avec le degr{\'e} et la  hauteur globale de la sous-vari{\'e}t{\'e} en question. 
Le lien est donn{\'e} par le {\it 
th{\'e}or{\`e}me des minimums successifs}~\cite[Thm.~5.2]{Zhang95}:
pour
une sous-vari{\'e}t{\'e} irr{\'e}ductible  $V \subset \P^N$ de dimension
$r$  
\begin{equation} \label{zhang}
\mu_1 (V)+ \cdots + \mu_{r+1} (V) 
\, \le \,  
\frac{h(V) }{\deg (V)}
\, \le \,  (r+1) \, \mu_1(V).   
\end{equation} 
On renvoie 
{\`a}~\cite[Thm.~3.1]{DaPh98}
pour une d{\'e}monstration {\'e}l{\'e}mentaire, 
bas{\'e}e sur 
les th{\'e}or{\`e}mes de
B{\'e}zout et de Hilbert-Samuel arithm{\'e}tiques.  

\smallskip 

Tout d'abord, il est naturel 
de se demander si ces estimations sont pr{\'e}cises pour
des cas particuliers. 
Le th{\'e}or{\`e}me~\ref{1} montre que d{\'e}j{\`a} pour les espaces projectifs, 
aucune de ces  estimations  n'est
exacte: 
on a 
$$
\mu_1 (\P^n)+ \cdots + \mu_{n+1} (\P^n) \, <   \, 
{h(\P^n) }
\, <  \, (n+1) \, \mu_1(\P^n),   
$$
car \ 
$ \displaystyle{ h(\P^n) =
\frac{n+1}{2}  \, \sum_{j=2}^{n+1} \frac{1}{j}  } 
$
\cite[Lem.~3.3.1]{BoGiSo94}, et d'apr{\`e}s le th{\'e}or{\`e}me~\ref{1} 
\vspace{-3mm} 
$$ 
\mu_1 (\P^n)+ \cdots + \mu_{n+1} (\P^n) =
\frac{1}{2} \, \log (n+1)
+ \cdots +  \frac{1}{2} \, \log (1) 
 \quad , \quad 
(n+1) \, \mu_1 (\P^n)
= \frac{n+1}{2} \, \log (n+1). 
 $$

Le r{\'e}sultat suivant est une cons{\'e}quence directe du th{\'e}or{\`e}me~\ref{1} et
des
estimations~(\ref{zhang}): 

\begin{cor}{(du th{\'e}or{\`e}me~\ref{1})}  \label{altura}
Soit $\cA =\{ a_0, \dots, a_N \} \subset \Z^n$ un ensemble fini 
de dimension $r$, alors 
$$
\frac{1}{2}  
 \, \Big( \log N_\cA(r) + \cdots +  
\log N_\cA(0) \Big) 
\, \Vol (\cA) 
\, \le \,
h(X_\cA)
\, \le \,
\frac{(r+1)}{2} \, \log(N+1) \, \Vol (\cA).
$$
\end{cor}

Soit $\cA= \{a_0, \dots, a_N \} \subset \Z^n$ un ensemble fini et 
$\alpha \in (\Qbarra^\times)^{N+1}$ un vecteur {\`a} coordonn{\'e}es
alg{\'e}briques, 
et consid{\'e}rons l'application monomiale
$$
\varphi_{\cA, \alpha} : \T^n \to \P^N 
\quad \quad , \quad \quad 
t \mapsto (\alpha_0 \, t^{a_0} : \cdots : \alpha_{N} \, t^{a_N}).
$$
La {\it vari{\'e}t{\'e} monomiale} $X_{\cA,\alpha} \subset \P^N$ associ{\'e}e est 
d{\'e}finie comme  {\'e}tant la cl{\^o}ture de Zariski de l'image de cette application,
c'est-{\`a}-dire
$ X_{\cA,\alpha} :=  \overline{\varphi_{\cA,\alpha}(\T^n)}$. 
Les vari{\'e}t{\'e}s toriques projectives sont des cas particuliers 
des vari{\'e}t{\'e}s monomiales (correspondant au cas $\alpha=(1,\dots,1)$)
et quelques unes des estimations pr{\'e}c{\'e}dentes 
peuvent s'{\'e}tendre sans difficult{\'e}  
{\`a} cette situation plus g{\'e}n{\'e}rale.

\smallskip 

Notons que 
$X_\arith$ et $X_\cA$ 
sont lin{\'e}airement isomorphes par l'application diagonale
$$
\P^N \to \P^N  
\enspace, \quad \quad 
(x_0 : \cdots : x_N ) \mapsto
({\alpha_{0}^{-1}} \, x_{0} : \cdots : 
{\alpha_{N}^{-1}} \, x_{N} ),
$$
et donc elles jouissent des m{\^e}mes propri{\'e}t{\'e}s {g{\'e}om{\'e}triques}; 
en particulier $ \dim(X_\arith) = \dim(\cA)$
et $\deg(X_\arith) = \Vol(\cA)$. 
Pour  une face 
\,$P \in \Faces(Q_\cA ) $\,  du polytope 
$Q_\cA \subset \R^n$, 
on notera  
$$
\alpha(P):= ( \alpha_i \, ; \ a_i \in P ) \in
(\Qbarra^\times)^{N(P)+1}.
$$ 

\begin{prop} \label{monomial} 
Soit $\cA= \{a_0, \dots, a_N \} \subset \Z^n$ un ensemble fini
de dimension $r$ et $\alpha \in (\Qbarra^\times)^{N+1}$, alors

\medskip 

\begin{itemize} 

\item[$\bullet$] 

\ $
\mu_i(X_\arith) \, \le \, 
 \min 
\,\Big\{ \,h(\alpha(P)) \ ; \ P \in \Faces(Q_\cA) , \ \dim (P) = r-i+1 
\,\Big\} 
$ \ 
pour $i=1, \dots, r+1$, et   

\medskip

\item[$\bullet$]
\ $ \displaystyle{ h(X_\arith) \, \le \,   \frac{r+1}{2}  \, h(\alpha)  \, \Vol
  (\cA) }$. 

\end{itemize} 
\end{prop} 

La d{\'e}monstration suit exactement la d{\'e}marche de  
celles  du th{\'e}or{\`e}me~\ref{1} et du
corollaire~\ref{altura}. 
Remarquons 
que la majoration pour les minimums successifs 
n'est plus une 
{\'e}galit{\'e} dans le cas g{\'e}n{\'e}ral
(consid{\'e}rer l'exemple \,$\cA:= \{ 0,1\} \in \Z$\, et \,$\alpha:= (1,2) \in \Q^2$).

\medskip

Plus g{\'e}n{\'e}ralement, cette m{\'e}thode nous permet d'estimer le 
comportement
du minimum essentiel et de la hauteur des vari{\'e}t{\'e}s
 par rapport {\`a} des morphismes.
Soit \,$\varphi : \P^N \to \P^M$\, une application rationnelle d{\'e}finie
par des formes 
\,$\varphi_0, \dots, \varphi_M  \in \Qbarra[x_0, \dots, x_N]$\,
de degr{\'e} $d$.
On d{\'e}finit son {\it degr{\'e}}  par \,$\deg (\varphi) := d$, et
sa {\em hauteur} \,$h_\W(\varphi) $ par la formule
\begin{eqnarray*} 
h_\W(\varphi)  & :=   & 
\sum_{v \in M_K^\infty} 
\lambda_v(K)  \, \log \sqrt{ 
||\varphi_0||_{\W, v}^2 + \cdots 
+ ||\varphi_M||_{\W,  v}^2 } 
\hspace{2mm} \\[0mm] 
&&+ \hspace{-2mm}
 \sum_{v \in  M_K \setminus 
M_K^\infty} \lambda_v(K)  \, \log \max \Big\{ |\varphi_0|_v, 
\dots, | \varphi_M|_v \Big\},  
\end{eqnarray*}
o{\`u}  $K $ est un corps de nombres 
contenant les coefficients de
$\varphi$;
pour  
$ v\in M_K^\infty$
on d{\'e}signe par 
\,$|| \cdot ||_{\W, v} $\, 
la m{\'e}trique de Weyl 
relative {\`a} la place $v$, {\it voir} le paragraphe~\ref{points}.

\begin{prop}
Soit \  $\varphi : \P^N \to \P^M$ \ une application rationnelle, 
\ $V \subset \P^N$ \ une sous-vari{\'e}t{\'e} irr{\'e}ductible de
dimension $r$ et  \ $Z:= \overline{\varphi(V)} \subset \P^M$ \  
l'image de $V$  par $\varphi$, alors
$$
\mu^\ess (Z ) 
\le   h_\W(\varphi) +  \deg ( \varphi) \,  \mu^\ess (V) 
\quad \quad , \quad \quad   
\frac{ h(Z)}{ \deg (Z)}  
 \le 
(r+1) \, \bigg( h_\W(\varphi) +  \deg (\varphi) \, \frac{h(V)}{\deg (V)}
\bigg).
$$

\end{prop}

\begin{demo}
Soit $\xi \in \Qbarra^{N+1}$. 
Pour   $v \in M_K^\infty$ 
$$
\log || \varphi(\xi) ||_{2, v} 
\, \le \,  
\log \sqrt{ ||\varphi_0||_{\W, v}^2 + \cdots + ||\varphi_M||_{\W, v}^2 }
+  d \, \log||\xi||_{2, v}
$$
comme cons{\'e}quence de l'in{\'e}galit{\'e}~(\ref{weyl}), 
tandis que 
pour $v \in M_K \setminus M_K^\infty$ on a
$$
\log \max \Big\{  |\varphi_0 (\xi) | , \dots , | \varphi_M(\xi) | 
\Big\} 
\, \le \,  
\log \max \Big\{ |\varphi_0|_v, 
\dots, | \varphi_M|_v \Big\}
+ 
d \,  \log \max \Big\{  |\xi_0|_v , \dots , |\xi_N|_v | \Big\} 
$$
gr{\^a}ce {\`a} 
l'in{\'e}galit{\'e} 
ultram{\'e}trique. 
On en obtient
\, $ h \Big(\varphi(\xi)\Big) \le h_\W(\varphi)+  d \, h(\xi) $ \, par
sommation sur $M_K$.

\smallskip

Maintenant soit  $\varepsilon > 0$  quelconque,  et soit \,$Y  
\subset Z$\, une sous-vari{\'e}t{\'e} propre tel que 
 \ $ \mu^\ess (Z) \le \mu^\abs(Z \setminus Y) + \varepsilon$.
Consid{\'e}rons alors  la cl{\^o}ture de Zariski   
\,$X:=  \overline{\varphi^{-1} (Y)} \subset \P^N $, qui est  
une sous-vari{\'e}t{\'e} propre de $V$ tel que \,$\varphi(X) = Y$. 
On a 
$$
 \mu^\abs(Z \setminus Y) \le 
h_\W(\varphi) 
+  d \, \mu^\abs(V \setminus X )  
\, \le \,  
h_\W(\varphi) + 
d \, \mu^\ess (V),
$$
ce qui d{\'e}montre la majoration pour le minimum essentiel. 
L'estimation pour le quotient  \,$h(Z) / \deg (Z)$\, suit directement 
de cette majoration et des estimations~(\ref{zhang}). 
\end{demo}

\medskip

Comme une autre  cons{\'e}quence du th{\'e}or{\`e}me des minimums successifs, on d{\'e}duit 
la minoration suivante pour la
hauteur d'une vari{\'e}t{\'e} projective  non
contenue dans la r{\'e}union des hyperplans coordonn{\'e}s: 

\begin{cor}{(de la proposition~\ref{bogomolov}(a))}   \label{inferieure}
Soit $V \subset \P^N$ une sous-vari{\'e}t{\'e} irr{\'e}ductible tel que 
 \ $ V \not\subset Z(x_0 \cdots x_N)$, alors 
$$
h(V) \ge \frac{1}{2} \, \log (N+1) \, \deg ( V). 
$$
\end{cor} 

Posons $r:= \dim(V)$. 
Lorsque la codimension de $V$ est grande, et plus pr{\'e}cis{\'e}ment lorsque 
$\displaystyle \frac{1}{2} \, \log (N+1) <  h(\P^r) = \frac{r}{2} \,
(\log(r+1) + O(1))$,  
ce r{\'e}sultat  am{\'e}liore la minoration 
\ 
$ h(V) \ge h(\P^r) \, \deg (V) $ ($r:=\dim (V)$) 
due {\`a} J.-B. Bost, H. Gillet et C. Soul{\'e}~\cite[Prop.~4.1.2(i) et Thm.~5.2.3]{BoGiSo94}.
Notons que cette minoration est {\'e}quivalente {\`a} la 
 positivit{\'e} de la hauteur des sous-vari{\'e}t{\'e}s  consid{\'e}r{\'e}e dans cette
 r{\'e}f{\'e}rence.

\medskip

Consid{\'e}rons maintenant  la majoration pour la taille des coefficients 
du r{\'e}sultant creux. 

Soit $\cA:= \{ a_0, \dots, a_N \} \subset \Z^n$ un ensemble fini tel
que 
\,$L_\cA = \Z^n$.
Pour chaque \,$i=0, \dots, n$,  on introduit  
un groupe de $N+1$ variables
 \,$U_i=\{ U_{i\,0}, \dots, U_{i\,n}\}$\, et 
posons
$$
F_i:= \sum_{j=0}^N U_{i \, j } \, x^{a_j} \  \in \Q[U_i] [x_1^{\pm 1} ,
\dots, x_n^{\pm 1} ] 
$$
pour le polyn{\^o}me de Laurent g{\'e}n{\'e}rique de support $\cA$. 
Soit 
$$
\Omega_\cA:= \Big\{ 
 (\nu_0, \dots, \nu_n; \, \xi) \in 
(\P^{N})^{n+1} \times \T^n \ ;   \ F_i (\nu_i, \xi)=0,  \ i = 0,
\dots, n 
\Big\} \subset (\P^N)^{n+1} \times \T^n
$$
la vari{\'e}t{\'e} d'incidence de \,$F_0, \dots, F_n$\, sur $\T^n$,
et soit  \,$  \pi : (\P^{N})^{n+1}
\times \T^n \to (\P^{N})^{n+1}  $\,
la projection canonique. 
Alors 
\,$ \overline{\pi (\Omega_\cA )} \ \subset(\P^{N})^{n+1} $\,
est une sous-vari{\'e}t{\'e} irr{\'e}ductible 
de codimension 1, 
et le
 {\it $\cA$-r{\'e}sultant} (ou r{\'e}sultant creux) $\Res_\cA$
est d{\'e}fini comme {\'e}tant l'unique ({\`a} un signe pr{\`e}s) 
polyn{\^o}me irr{\'e}ductible  
d{\'e}finissant cette  
hypersurface \  \cite[Ch. 8,
Prop.-Defn. 1.1]{GeKaZe94}. 

En fait, 
le $\cA$-r{\'e}sultant co{\"\i}ncide avec la forme de Chow de
la vari{\'e}t{\'e} torique projective
 $X_\cA$ \ 
\cite [Ch. 8, Prop. 2.1]{GeKaZe94}. 
On peut estimer  de fa{\c c}on routini{\`e}re 
la taille maximale $h_\mmax (\Res_\cA)$ de ses coefficients, 
{\`a} partir de la majoration pour $h(X_\cA)$.

\begin{Demo}{D{\'e}monstration du corollaire~\ref{res}.--} 
Le $\cA$-r{\'e}sultant est , par d{\'e}finition,   
un  polyn{\^o}me primitif (c'est-{\`a}-dire ses coefficients sont premiers entre
 eux)
et donc de hauteur locale nulle
pour toutes les  places ultram{\'e}triques. 
Ainsi
$$
h(X_\cA) = m\Big(\Res_\cA \, ; \, S_{N+1}^{n+1} \Big) + (n+1)  \, 
\bigg(\sum_{i=1}^N \frac{1}{i}\bigg)  \, \Vol(\cA), 
$$
o{\`u}  \,$ m\Big(\Res_\cA \, ; \, S_{N+1}^{n+1} \Big) $\, d{\'e}signe la
$S_{N+1}^{n+1} $-mesure du polyn{\^o}me $\Res_\cA$, {\it voir} 
le paragraphe~\ref{points}. 
Donc 
$$\begin{array}{rcl} 
 h_\mmax(\Res_\cA) 
& \le & 
m(\Res_\cA  ) 
+ (n+1) \, \log (N+1) \, \Vol(\cA) \\[3mm] 
& \le & 
h(X_\cA) 
+ (n+1) \, \log (N+1) \, \Vol(\cA) \\[3mm] 
& \le & \displaystyle \frac{3}{2} \, (n+1) \, \log (N+1) \, \Vol(\cA), 
\end{array}$$
o{\`u}  $m(\Res_\cA)$ d{\'e}signe la mesure de Mahler de $\Res_\cA$. 
La premi{\`e}re in{\'e}galit{\'e} est  cons{\'e}quence de~\cite[Lem.~1.1]{KrPaSo01}
en regardant 
$\Res_\cA$ comme un polyn{\^o}me en $n+1$ groupes de $N+1$
variables chacun, homog{\`e}ne de degr{\'e} $\Vol(\cA)$ dans chaque groupe. 
La deuxi{\`e}me et la troisi{\`e}me in{\'e}galit{\'e}s sont des cons{\'e}quences 
de~\cite[Ineq.~(1.2)]{KrPaSo01} et du 
corollaire~\ref{altura}, respectivement. 
\end{Demo}

% ----------- Seccion 3 -----------------------------------------------

\typeout{Koushnirenko} 

\section{Un analogue arithm{\'e}tique du 
th{\'e}or{\`e}me de Koushnirenko}

\label{koushnirenko}

\setcounter{equation}{0}

Dans ce paragraphe on 
d{\'e}montre  un analogue arithm{\'e}tique du 
th{\'e}or{\`e}me de Koushnirenko pour la hauteur des solutions  
d'un syst{\`e}me d'{\'e}quations polynomiales. 
On pr{\'e}sente aussi des variantes de ce r{\'e}sultat pour 
le cas des
intersections impropres.

\bigskip

Soit \,$\cA := \{ a_0, \dots, a_N \}  \subset \Z^n$\, un ensemble
fini. 
Pour une  
sous-vari{\'e}t{\'e} irr{\'e}ductible 
\,$ V \subset \T^n  $, 
on   
consid{\`e}re son   {\em $\cA$-degr{\'e}}  et sa  {\em $\cA$-hauteur}, 
respectivement d{\'e}finis par  
$$
\deg_{\cA}  (V) := \deg (\varphi_{\cA} (V)) 
\quad \quad , \quad \quad
h_\cA (V) := h(\varphi_{\cA} (V)).   
$$

On {\'e}tend par lin{\'e}arit{\'e} ces d{\'e}finitions 
au groupe des cycles \,$Z(\T^n)$\, du tore
$\T^n$. 
Pour une sous-vari{\'e}t{\'e} quelconque 
$V \subset \T^n$, on d{\'e}finit son   
 {\em $\cA$-degr{\'e}}  et sa {\em $\cA$-hauteur}  comme
ceux de son cycle associ{\'e}.

\medskip

Dans la suite on d{\'e}finit  un 
produit d'intersection entre cycles et diviseurs de  $\T^n$. 
Tout d'abord, remarquons que tout diviseur de Cartier 
\,$D \in  \Div(\T^n) $\, 
est principal car l'anneau  
$ \Qbarra[x_1^{\pm 1} , \dots, x_n^{\pm 1}]$ est factoriel. 
Posons alors 
\,$f_D \in \Qbarra[x_1^{\pm 1} , \dots, x_n^{\pm 1}]$\, 
pour le polyn{\^o}me de Laurent (unique {\`a} un facteur scalaire pr{\`e}s) 
d{\'e}finissant $D$. 

Soit \,$V \subset \T^n$\, 
une sous-vari{\'e}t{\'e} irr{\'e}ductible   non
contenue dans le support 
 $|D |  \subset \T^n$ du diviseur $D$. 
Pour chaque composante irr{\'e}ductible \,$C \in \Irr(V \cap  |D|) $\, 
de l'intersection ensembliste  $V \cap |D|$
on consid{\'e}rera 
la multiplicit{\'e} d'intersection classique d{\'e}finie par 
la {longueur}
$$
\ell(V, D \, ; \, C) := \lg \left( \Qbarra[V] 
/ (f_D)\right)_{ I(C)}, 
$$
o{\`u}  \,$I(C) \subset
R: = \Qbarra[x_1^{\pm 1} , \dots, x_n^{\pm 1}]$\,  
d{\'e}signe l'id{\'e}al premier de d{\'e}finition de 
$C$ 
et \,$\lg $\, la longueur du  
$R_{I(C)}$-module $ \left( \Qbarra[V] 
/ (f_D)\right)_{I(C)} $. 
Pour une sous-vari{\'e}t{\'e} irr{\'e}ductible    
  \,$V \subset \T^n$\, et un diviseur   
\,$D \in  \Div(\T^n) $\, quelconques on d{\'e}finit le produit  
$$
[V] \cdot D := \left\{ 
\begin{array}{ll} 
 \sum_{C \in \irr(V \cap  |D|) } 
\, \ell(V,D \, ; \, C) \, [C] \quad \quad \quad  
\quad & \mbox{ si } V  \not\subset |D|, \\[3mm] 
[V] & \mbox{ sinon,} 
\end{array} \right. 
$$
et on l'{\'e}tend
par lin{\'e}arit{\'e} en un accouplement
$$
Z(\T^n) \times \Div(\T^n) \to Z(\T^n)
\quad \quad , \quad \quad
(Z, D) \mapsto Z \cdot D.   
$$
Pour plusieurs diviseurs \,$D_1, \dots, D_s$\,  
on pose 
\ $ Z \cdot D_1 \cdot D_2 \cdots   D_s  := 
\Big( (Z \cdot D_1) \cdot D_2 \Big) \cdots  D_s$; 
notons que ce cycle {\it d{\'e}pend} de l'ordre des diviseurs choisi. 
  
\smallskip

\begin{defn} 
Soit $Z \in Z(\T^n)$ un cycle effectif et
$D_1, \dots,  D_s \in \Div(\T^n)$ des diviseurs   
effectifs, alors 
$$
Z \cdot D_1 \cdots D_s  = \sum_C 
m(C) \, [C],  
$$
o{\`u} $C$ parcours l'ensemble des sous-vari{\'e}t{\'e}s irr{\'e}ductibles de  
l'intersection ensembliste $ |Z| \cap |D_1| \cap \cdots \cap |D_s|  $ 
et $m(C) \in \N$
avec  $m(C) =0$ sauf pour une nombre fini des $C$. 
Cet entier 
$$ m(Z, D_1, \dots, D_s
\, ; \, C) := m(C) 
$$ 
est par d{\'e}finition   
la {\em multiplicit{\'e} d'intersection}  de $Z$ avec $D_1, \dots, D_s$
le long de $C$. 
\end{defn}

On v{\'e}rifie ais{\'e}ment \,$m(C) \ge 1$\, lorsque $C$ est une
composante isol{\'e}e de $ |Z| \cap |D_1| \cap \cdots
  \cap |D_s| $, et donc 
le degr{\'e} et la hauteur de l'intersection ensembliste
sont major{\'e}s par ceux 
 du cycle
intersection: 
\begin{eqnarray*} 
\deg_\cA (|Z| \cap |D_1| \cap
\cdots \cap |D_s|)
 &  \le &   \deg_\cA (Z \cdot D_1 \cdots D_s)  , \\[2mm] 
 h_\cA (|Z| \cap |D_1| \cap
\cdots \cap |D_s|)
 & \le & 
h_\cA (Z \cdot D_1 \cdots D_s).
\end{eqnarray*}

Par la suite on 
montre que pour une composante {\it propre}
$C$  de
l'intersection d'une famille de diviseurs de  $\T^n$, cette multiplicit{\'e}
co{\"\i}ncide
avec la longueur; 
on d{\'e}montre aussi la multilin{\'e}arit{\'e} de $m(C)$. 

\begin{lem} \label{longitud} 

Soient 
 \,$ D_1, \dots, D_s \in \Div(\T^n)$ et 
soit \,$C \subset \T^n$\, 
une composante irr{\'e}ductible de dimension $n-s$ de 
$ |D_1| \cap \cdots \cap |D_s|$, alors 

\medskip 

\begin{itemize} 

\item[(a)]  
\ $ 
m(\T^n, D_1, \dots,  D_s \, ; \, C) = 
\lg \left( \Qbarra[x_1^{\pm 1} , \dots, x_n^{\pm
 1}]  /(f_{D_1} , \dots, f_{D_s})\right)_{I(C)} $, et 

\smallskip 

\item[(b)] l'application \,$\Z^s \to \Z$, 
\,$ (k_1, \dots, k_s)  \mapsto 
m(\T^n, k_1 \, D_1, \dots,  k_s \, D_s \, ; \, C) $\, est
multilin{\'e}aire. 

\end{itemize} 

\end{lem}

\begin{demo}
Consid{\'e}rons d'abord la partie (a), qu'on d{\'e}montrera par 
r{\'e}currence en $s$. 
Le cas $s=0$ {\'e}tant trivial, on consid{\`e}re 
le cas 
$s \ge 1$ en 
supposant que l'{\'e}nonc{\'e} est valable 
pour $s-1$. 

\medskip 

Soit 
\ $R: = \Qbarra[x_1^{\pm 1} , \dots, x_n^{\pm 1}]$, 
 \ $\ga:= (f_{D_1}, \dots, f_{D_{s-1}}) \subset R$, 
\ $ \gp : = I(C)  \subset R$, 
et 
\ $ 
A:= ( R/ \ga )_{\gp}    $. 
Autrement dit, 
\,$\Big(A, (\gp)  \Big)$\,  est l'anneau local $  \cO_{V, C}$ de $V$ 
le long de
$C$. 
On a $\dim(\ga)= n-s+1$ et $\dim(\gp)= n-s= \dim(\ga)-1$, ce
qui entra{\^\i}ne que $A$ est un anneau de dimension 1. 
En outre, le polyn{\^o}me de Laurent $f:= f_{D_s} \in R $ appartient {\`a} $\gp$ et 
c'est un non-diviseur de z{\'e}ro de $A$. 
On en d{\'e}duit que 
$A$ est Cohen-Macaulay car en dimension 1, un anneau est Cohen-Macaulay si et seulement si
il contient un non diviseur de z{\'e}ro.  

Aussi on a  \,$\dim
\Big(A/(f)\Big)  =0$ car 
$f$ est un non-diviseur de z{\'e}ro;  
on  consid{\`e}re alors  la multiplicit{\'e} de Samuel
\,$e\Big((f) , A\Big)$\, de  $f$ dans $A$~\cite[Sec. 1.2]{FlOCVo99}.

Soit \,$\gq \in \ass(A) $\, un id{\'e}al premier associ{\'e} de  $A$  
quelconque. 
Similairement on d{\'e}montre que \,$\Big(A / \gq, (\gp)\Big) $\, est un anneau local
Cohen-Macaulay de dimension 1; en fait c'est un {\it domaine}. 
{\`A} nouveau \,$\dim \Big(A/ (\gq +f) \Big) = 0$\, et on 
consid{\`e}re aussi la 
multiplicit{\'e} $e\Big( (f) , A/\gq\Big)$ de $f$ dans $A/\gq$.

\smallskip 

On v{\'e}rifie que $f$ est un $G$-param{\`e}tre
pour $\Big( A, (\gp) \Big) $ et pour  $\Big( A/\gq, (\gp) \Big)$, 
{\it voir}~\cite[Defn. 1.2.10]{FlOCVo99}, 
et donc  
$$
e\Big( (f)  , A\Big)  = \lg \Big(  A/(f) \Big)   
\quad \quad , \quad \quad 
e \Big( (f) ,  A/ \gq \Big) = \lg \Big(  A/(\gq + f)\Big)
$$
gr{\^a}ce {\`a}~\cite[Cor.~1.2.13]{FlOCVo99}. 
La 
formule d'associativit{\'e}~\cite[Thm.~1.2.8]{FlOCVo99} \ 
entra{\^\i}ne 
$$
 \lg \Big(  A/(f) \Big)   
= e \Big( (f)  , A \Big)
=  \sum_{\gq\in \Ass(A) } \lg (A_\gq)    
\, e \Big( (f) ,  A/ \gq  \Big) 
=  \sum_{\gq\in \Ass(A) } \lg (A_\gq)    
\, \lg \Big(  A/(\gq +  f)\Big),
$$
et on en conclut
\begin{eqnarray*} 
m(C) 
& = & 
\sum_{\gq\in \Ass(A) } 
m \Big(\T^n, D_1, \dots,  D_{s-1}  \, ; \, Z(\gq) \Big) \ 
\ell \Big( Z(\gq) , D_s \, ; \, C \Big) 
\\[2mm] 
& = & 
  \sum_{\gq\in \Ass(A)}  \lg (A_\gq)    
\, \lg \Big(  A/(\gq +  f)\Big)  \\[2mm] 
& = &    
\lg \Big(  A/(f) \Big),  
\end{eqnarray*} 
car 
\ $m \Big(\T^n, D_1, \dots , D_{s-1} \, ; \, Z(\gq) \Big)=
 \lg (A_\gq)$ \ 
par l'hypoth{\`e}se de r{\'e}currence. 

\medskip

D'apr{\`e}s ce qu'on vient de voir, la partie (b) se ram{\`e}ne {\`a} 
v{\'e}rifier  \ $ e \Big( (f^k)  , A \Big) = k\, e \Big( (f)  , A
\Big)$
pour tout $k \in \N^\times$. 

Ceci est une cons{\'e}quence directe des  d{\'e}finitions: 
soit \,$P_A^{(f)} $\, et \,$  P_A^{(f^k)}$\, 
les fonctions de Hilbert-Samuel de
$ \Big((f) , A\Big) $ et de  $ \Big( (f^k)  , A\Big)$ respectivement, 
{\it voir}~\cite[\S~1.2]{FlOCVo99}. 
Alors
$$ 
P_A^{(f^k)} (t) = \lg \Big(A/ (f^k)^t \Big)
= \lg  \Big(A/ (f^{k\,t} ) \Big)    
= P_A^{(f)} (k \, t) 
$$
et donc 
 \ $  
e\Big( (f^k)  , A\Big) \, t + O(1)= P_A^{(f^k)} (t)  = P_A^{(f)} (k \, t)
=  e\Big( (f^k)  , A\Big) \, k \, t + O(1)$ pour $t \gg0$, 
ce qui {\'e}tabli 
l'{\'e}galit{\'e} cherch{\'e}e. 
\end{demo}

\smallskip

Pour un  polyn{\^o}me de Laurent 
$f \in \Qbarra [x_1^{\pm 1} , \dots, x_n^{\pm 1} ]$, on note 
\,$ h_2 (f)$\, la hauteur
associ{\'e}e {\`a} la norme $\ell^2$, {\it voir} le
paragraphe~\ref{points}.

\begin{lem}  \label{BK-dot} 

Soit  $ \cA \subset \Z^n $ un ensemble fini tel que
 $L_\cA = \Z^n$,  
\,$f_1, \dots,  f_s \in
\Qbarra[x_1^{\pm 1} , \dots, x_n^{\pm 1}]$\, des polyn{\^o}mes de 
Laurent 
tels que \,$\Supp(f_i) \subset \cA$\, pour $i=1, \dots, s$ et 
\,$Z \in Z(\T^n)$\, un cycle effectif, alors 
\begin{eqnarray*} 
\deg_\cA \Big(Z \cdot \div(f_1)  \cdots \div(f_s) \Big) 
& \le & 
  \deg_\cA ( Z) , \\
h_\cA \Big(Z \cdot \div(f_1) \cdots  \div(f_s) \Big)  
& \le &  
h_\cA(Z) +  \deg_\cA(Z) \, \sum_{i=1}^s h_2(f_i). 
\end{eqnarray*} 

\end{lem} 

\begin{demo}
Il suffit de d{\'e}montrer le cas $s:=1$ de l'{\'e}nonc{\'e}, 
car le cas g{\'e}n{\'e}ral s'en suit par it{\'e}ration. 
En plus, par lin{\'e}arit{\'e}
on peut se ramener au cas o{\`u} \,$Z =[V]$, 
o{\`u} $V \subset \T^n$   est une sous-vari{\'e}t{\'e}  
irr{\'e}ductible.
Le cas  $V \subset Z(f_1) $ {\'e}tant  {\'e}vident,  on 
 suppose  sans perte de g{\'e}n{\'e}ralit{\'e}  $V $ non contenue dedans $Z(f_1)$. 

\smallskip 

Soit \,$W:= {\varphi_\cA (V) } \subset \P^N $\, avec $N+1 = \Card(\cA)$, 
et notons
\,$ 
\ell:=  (\varphi_\cA^*)^{-1}(f)  \in \Q [y_0,
\dots, y_N ]$\, la forme lin{\'e}aire correspondant 
au polyn{\^o}me de Laurent $f:= f_1$.

L'hypoth{\`e}se $ L_\cA = \Z^n$ 
 {\'e}quivaut {\`a} ce que 
\,$ \varphi_\cA : \T^n \to X^\circ_\cA$\, soit un isomorphisme, 
ce qui implique que la restriction  
\,$
\varphi_\cA : V \to W $\, 
est aussi un isomorphisme. 
On en d{\'e}duit   
l'{\'e}galit{\'e} des cycles 
$$
(\varphi_\cA)_* \Big( V \cdot \div(f) \Big):= 
\sum_{C \in \irr\big(V \cap Z(f)\big)} \ell(V, \div(f) \, ; \, C) \,
 [ \varphi_\cA(C) ] 
= W \cdot \div(\ell)  \ \in Z(\T^n), 
$$ 
et donc 
\ $
\deg_\cA\Big (V \cdot \div(f)\Big) = 
\deg  \Big( W \cdot \div(\ell)\Big) $ \ et 
\ $ 
h_\cA\Big (V \cdot \div(f)\Big) = h\Big ( W \cdot \div(\ell)\Big) $. 

\smallskip

La cl{\^o}ture de Zariski \,$\overline{W} \subset \P^N$\, 
n'est pas contenue dans $Z(\ell)$ et donc  
$$
\deg  \Big( \overline{W} \cdot \div(\ell)  \Big)  =  \deg ( \overline{W})
\quad \quad , \quad \quad 
h  \Big(\overline{W} \cdot \div(\ell)  \Big) 
=    h(\overline{W} ) +  \deg (\overline{W}) \, h_{\overline{W}} (\ell) , 
$$
par le 
th{\'e}or{\`e}me de B{\'e}zout arithm{\'e}tique dans la version  
de~\cite[Prop.~4]{Philippon919495}.  
Ici 
\,$h_{\overline{W}} (\ell)$\, d{\'e}signe la hauteur de 
$\ell$ relative {\`a}  la sous-vari{\'e}t{\'e} $\overline{W}$,
on renvoie {\`a}~\cite[p.~355]{Philippon919495} pour sa d{\'e}finition pr{\'e}cise. 
Cette hauteur est born{\'e}e par  
 \,$
h_{\overline{W}} (\ell)  \le h_2(\ell) $, ce qui est {\'e}tabli au cours de la  
d{\'e}monstration de~\cite[Prop.~4]{Philippon919495}.  

\smallskip

En outre, le cycle 
\,$  W \cdot \div(\ell) $\, est la restriction de 
\,$ \overline{W} \cdot \div(\ell) $\, {\`a} l'ouvert
$ \T^N \hookrightarrow  \P^N$, et donc 
$$
 \deg  \Big ( W \cdot \div(\ell) \Big) 
\, \le \, 
\deg  \Big( \overline{W} \cdot \div(\ell)  \Big)  \, = \,   
\deg ( \overline{W}) \, = \,  \deg_\cA (V),  
$$
et similairement 
$$
h   \Big( W \cdot \div(\ell) \Big)  
\, \le \,  
h  \Big(\overline{W} \cdot \div(\ell)  \Big) 
\, \le \,  h(\overline{W} ) + \deg (\overline{W})  \, h_2 (\ell)  
\, =  \, h_\cA(V) + \deg_\cA(V)  \, h_2(f).  
$$
Ainsi on a {\'e}tabli  \ $ \deg_\cA \Big(V \cdot \div(f) \Big ) \, \le \, 
\deg_\cA (V)$ \ 
et \ $ h_\cA \Big(V \cdot \div(f) \Big) \, \le \, h_\cA(V) 
+ \deg_\cA (V) \, h_2(f) $.
\end{demo}

\smallskip

Soit \,$Q \subset \R^n$\, un polytope rationnel et 
\,$\cA_Q:= Q \cap \Z^n $\, 
l'ensemble de ses points entiers. 

La  $Q$-hauteur de Weil d'un point  $\xi \in \P^N$ 
se compare avec sa 
$\cA_Q$-hauteur projective: on a 
\begin{equation*} \label{maj} 
\widehat{h}_Q(\xi) \le h_\cA(\xi) \le \widehat{h}_Q(\xi) 
+ \frac{1}{2} \, \log (\Card(\cA)). 
\end{equation*}
L'application $Q \mapsto \widehat{h}_Q$ est additive 
par rapport {\`a} la somme de Minkowski des polytopes, 
et invariante par translations:

\begin{lem} \label{funtorial} 
\

\begin{itemize} 

\medskip 

\item[(a)] 
Soit  $P, Q
\subset \R^n$ des polytopes rationnels, alors 
\ $
\widehat{h}_{P+ Q} = \widehat{h}_P + \widehat{h}_Q$. 

\medskip 

\item[(b)] 
Soit  $ Q
\subset \R^n$  un polytope rationnel et
$b \in \Z^n$, alors \ $ \widehat{h}_{b+ Q} = \widehat{h}_Q$. 

\end{itemize} 

\end{lem} 

\begin{demo}
Soit $\cA \subset \Z^n$ un ensemble fini quelconque et 
$\widehat{h}_\cA$ la 
$\cA$-hauteur de Weil associ{\'e}e, {\it voir} le
paragraphe~\ref{points}. 
Soit 
\,$ c \in \Conv(\cA) \cap \Z^n$, et prenons des 
des r{\'e}els non-n{\'e}gatifs
\,$ \{ r_a \in \R_+ \, ; \, a \in \cA\} $\, tels que 
 \,$
\sum_{a \in \cA} r_a \, a = c$\, et \,$
\sum_{a \in \cA} r_a = 1$. 
Pour $\xi \in \T^n$ et  $v \in M_K$ on a
$$
|\xi^c|_v = \prod_{a \in \cA} |\xi^{a}  |_v^{r_a} 
\ \le  \  
\max \Big\{ |\xi^{a}  |_v \, ; \,  a \in \cA  \Big\}; 
$$
ainsi on a montr{\'e}  \,$\widehat{h}_\cA(\xi) = \widehat{h}_{\conv(\cA)}
(\xi)$. 

\smallskip 

Maintenant soit \,$\cB:= P \cap \Z^n$\, et \,$ \cA := Q \cap \Z^n$. 
La partie (a) se d{\'e}duit de l'identit{\'e} 
\,$ \Conv(\cB + \cA) = P+Q$\, due {\`a}~\cite[Ch.~7, Prop.~4.3]{CoLiOS98},
tandis que la partie (b) est cons{\'e}quence imm{\'e}diate  
de la partie (a) et de la
formule du produit. 
\end{demo} 

Consid{\'e}rons {\`a} pr{\'e}sent la hauteur \,$h_1(f)$\,
associ{\'e}e {\`a} la norme $\ell^1$, {\it voir} le paragraphe~\ref{points}. 
Cette hauteur se compare avec la hauteur $h_2$; on a  
$\displaystyle h_2(f) \le h_1(f) \le h_2(f) + \frac{1}{2} \,
\log\Big(\Card(\Supp(f))\Big)$.
En plus, elle  est sous-additive: 
pour \,$f, g \in \Qbarra [x_1^{\pm 1} , \dots, x_n^{\pm 1} ] $ on a  
$$ 
h_1(f \, g) \le h_1 (f) + h_1(g).
$$

\smallskip

Le th{\'e}or{\`e}me~\ref{2} est le cas $K:= \Q$ de l'{\'e}nonc{\'e} suivant:

\begin{thm} \label{2-general} 
Soient 
\,$f_1, \dots, f_n \in K[x_1^{\pm 1}, \dots, x_n^{\pm 1}] 
$\,
des polyn{\^o}mes de Laurent
{\`a} coefficients dans $K$, et posons 
\,$ Q:= \NP(f_1, \dots, f_n) \subset \R^n$. 

Soit 
\,$Z(f_1,\dots, f_n)_0 $\, l'ensemble des points isol{\'e}s de 
\,$Z(f_1, \dots, f_n) \subset \T^n$. 
Pour chaque point   $\xi$ dans cet ensemble on note
\ $\ell(\xi) := \dim_\Qbarra  
\left(
\Qbarra[x_1^{\pm 1} , \dots, x_n^{\pm 1} ]
/ (f_1, \dots, f_n) \right)_{I(\xi)}$\ la multiplicit{\'e}
d'intersection de
$f_1, \dots, f_n$ en  $\xi$.   
Alors
$$
\sum_{\xi \in Z(f_1,\dots, f_n)_0}  \ell(\xi) \, \widehat{h}_Q (\xi) 
\ \le \ 
n!\, \Vol_n (Q)  \, \sum_{i=1}^n \, h_1 (f_i). 
$$

\end{thm}

\begin{Demo}
Supposons pour le moment
\,$ L_{Q\cap \Z^n} = \Z^n$. 
Soit $ k \in \N^\times$, et posons 
   \,$\cA_k := (k \, Q)  \cap \Z^n$. 
Soit $ a \in Q \cap  \Z^n = \cA_1 $ 
un vecteur entier dans $Q$ quelconque, 
alors   \,$ (k-1) \, a + \cA_1 \subset \cA_k$\, et donc 
 \,$ L_{\cA_k} \supset L_{\cA_1} = \Z^n$\,  
par hypoth{\`e}se, 
ce qui entra{\^\i}ne \,$ L_{\cA_k} = \Z^n$.

Les polyn{\^o}mes de Laurent 
 $f_1^k, \dots, f_n^k$ sont support{\'e}s dans $\cA_k$; 
alors leur applique le lemme~\ref{BK-dot} et on 
 trouve 
$$
h_{\cA_k} \Big( \T^n \cdot \div(f_1^k) \cdots  \div(f_n^k) \Big)  
\ \le \ h_{\cA_k} (\T^n) +  \deg_{\cA_k} (\T^n) 
\, \sum_{i=1}^n h_2(f_i^k). 
$$

Soit \,$V := Z(f_1, \dots, f_n) \subset \T^n$ et posons 
\ $\ell_k(\xi) :=  \dim_\Qbarra  
\left(
\Qbarra[x_1^{\pm 1} , \dots, x_n^{\pm 1} ]
/ (f_1^k, \dots, f_n^k) \right)_{I(\xi)}$
\ 
pour chaque point isol{\'e} 
\,$\xi \in Z(f_1^k, \dots, f_n^k)_0  =V_0 $. 
La positivit{\'e} des 
multiplicit{\'e}s d'intersection et du 
lemme~\ref{longitud}(a) entra{\^\i}nent
$$
\sum_{ \xi \in V_0 } 
\ell_k(\xi) \, \widehat{h}_{k Q} (\xi) 
\ \le 
\ 
\sum_C m(C) \, h_{\cA_k} (C) = 
h_{\cA_k}  \Big(\T^n \cdot \div(f_1^k) \cdots  \div(f_n^k) \Big)  , 
$$ 
car 
 \,$\Conv(\cA_k) = k \, \Conv(\cA_1) = k \, Q$\ 
et 
\,$ \widehat{h}_{k Q} (\xi) \le h_{\cA_k} (\xi)  $. 
En outre,
 \,$\deg_{\cA_k} (\T) = \Vol(\cA_k) = n!\, \Vol_n(k \, Q) $\,
car $ L_{\cA_k} = \Z^n$ et on obtient 
$$
h_{\cA_k} (\T^n) +  \deg_{\cA_k} (\T^n) 
\, \sum_{i=1}^n h_2(f_i^k)  
\ \le \ 
n!\, \Vol_n (k \, Q)  \, \bigg( 
\frac{1}{2} \, (n+1) \, 
\log (\Card( \cA_k)  ) + 
\sum_{i=1}^n \, h_1 (f_i^k) 
\bigg)
$$
en utilisant aussi la majoration pour la hauteur de la vari{\'e}t{\'e} torique $X_{\cA_k} $ 
(corollaire~\ref{altura}) et l'estimation
\,$ h_2(f_i) \le h_1(f_i)$.

\smallskip 

Ensuite, 
on a  \,$ \ell_k(\xi) = k^n \, \ell (\xi) $\, par le 
lemme~\ref{longitud}(b) et  
\,$ \widehat{h}_{ k Q}  (\xi) = k\,  \widehat{h}_{ Q} (\xi) $\, 
gr{\^a}6ce {\`a} l'additivit{\'e}  de la hauteur $\widehat{h}$
\Big(lemme~\ref{funtorial}(a)\Big). 
Aussi on a   \,$\Vol_n(k \, Q) = k^n \, \Vol_n(Q)$\, et 
\,$ h_1 (f_i^k) \le k \,  h_1 (f_i)$. 

\smallskip 

Soient 
$ b \in \Z^n$ et 
$d \in \N^\times  $ tels que $Q \subset b + d \, S$, o{\`u} $S$ 
d{\'e}signe le
simplex
standard de $\R^n$. 
Alors \,$k \, Q \subset k \, b + k \, d \, S$\, et donc 
\,$
\log ( \Card( \cA_k) ) 
\le \log \Big( \Card( ( k \, b + k\, d \, S \, \cap \Z^n)) \Big) 
= \log  {k \, d+n \choose n } 
= O_k(\log k )$  (ici la notation  
$O_k$ r{\'e}f{\`e}re  {\`a} la d{\'e}pendance 
en $k$). 
On en obtient  
$$ 
\sum_{ \xi \in V_0 } 
\ell(\xi) \, 
\widehat{h}_{Q}  (\xi) 
\ \le \ 
n!\, \Vol_n (Q)  \, 
\sum_{i=1}^n \, h_1 (f_i)
+ 
O_k \bigg(\frac{\log k }{k}  \bigg), 
$$
d'o{\`u} on conclut  en faisant \ $ k \to \infty$. 

\medskip

Maintenant consid{\'e}rons le cas g{\'e}n{\'e}ral o{\`u} 
$\dim (Q) =n$. 
Soit   \,$ L_\cA$\,  le sous-module de $\Z^n$  
engendr{\'e} par les diff{\'e}rences des vecteurs dans
$\cA:= Q \cap \Z^n $, qui  
est un sous-r{\'e}seau de $\Z^n$, mais pas forcement 
{\em {\'e}gal}  {\`a} $\Z^n$. 
On montre dans la suite 
que ce cas se r{\'e}duit au cas pr{\'e}c{\'e}dent. 

\smallskip

Soient \,$c_1, \dots, c_n \in \N^\times$\, les facteurs {\'e}l{\'e}mentaires 
de $L_\cA$, et 
\,$v_1, \dots, v_n \in \Z^n$\, 
des vecteurs entiers formant une base de $\Z^n$ tels que 
\,$c_1 \, v_1 , \dots, c_n \, v_n$\, soit une base de  $L_\cA$. 

Soit \,$ \beta: \Z^n \to \Z^n$\, l'application lin{\'e}aire 
d{\'e}finie par $ v_i \mapsto c_i \,
v_i$, qui est un isomorphisme entre $\Z^n$ et $L_\cA$. 
Soit \ $\cB := \beta^{-1} (\cA) \subset \Z^n$\ 
et  \ $P:=\Conv (\cB)  \subset \R^n $. 
On a $L_\cB = \Z^n$ et donc 
  \,$
\Vol(\cB) = n!\, \Vol_n(P) =  n!\, \Vol_n(Q) / \gamma  
$,   
o{\`u} 
\,$ \gamma:= c_1 \cdots c_n$ 
est l'indice du r{\'e}seau  $L_\cA$.

\smallskip
 
Soit \,$\psi : \T^n \to \T^n $\, le morphisme  
d{\'e}fini par 
$ t 
\mapsto \Big(t^{\beta(e_1)}  , \dots, t^{\beta(e_n)} \Big)  $ 
o{\`u} les $e_j$ d{\'e}signent les vecteurs de la base standard de
$\R^n$, 
de telle sorte que \,$\varphi_\cA = \varphi_\cB \circ \psi$.

\smallskip 

Consid{\'e}rons l'inclusion d'anneaux 
\,$ \psi^*: B:= \Qbarra(\T^n) \hookrightarrow  
A:= \Qbarra(\T^n)$; 
on note que 
$$ 
\psi^*(B) = 
\Qbarra[x^{\pm \beta(e_1)} , \dots, x^{\pm \beta(e_n) } ] 
\ \subset \  
A= \Qbarra[x_1^{\pm 1} , \dots, x_n^{\pm 1}] 
$$
et donc $\psi^*$ est finie,
car
les  vecteurs $  {\beta(e_1)} , \dots,
{\beta(e_n) }$ sont lin{\'e}airement ind{\'e}pendants.

\smallskip 

Notons \,$K:= \psi^*(B)_{(0)}  \subset L:= A_{(0)} $\, 
les corps des fractions respectifs; 
on a  \ $\Card( \psi^{-1} (\xi)) = \gamma $\ pour tout 
$\xi \in \T^n$ et donc 
\,$ [L : K ] = \gamma $. 

\smallskip 

Soit \,$g_i := (\psi^*)^{-1}(f_i) \in B $\, pour $i=1, \dots, n$,
et  posons 
\,$\gothg:= (g_1, \dots, g_n) \subset B$\, et 
\,$\gf:= (f_1, \dots, f_n) = 
\Big( \psi^* (\gothg) \Big) \subset A $. 
Soit \,$W:= Z( \gothg) \subset \T^n$ la vari{\'e}t{\'e} d{\'e}finie par
$g_1, \dots, g_n$.  
Soit  $\eta \in W_0 $ un point isol{\'e}, 
$I(\eta) \subset B$ son id{\'e}al, et notons  
\,$\ell(\eta) = \lg  (B/\gothg)_{I(\eta)}  $\, 
la multiplicit{\'e} d'intersection 
de $g_1, \dots, g_n$ en $\eta$.

\smallskip

On remarque que  $B_{I(\eta)} $ est un anneau local, 
$A_{I(\eta)} $ est un $B_{I(\eta)} $-module fini de dimension $n$, 
et  
 \,$ \dim \Big( A_{I(\eta)} /  \gf \Big) = 0 $.
Par  la formule
d'associativit{\'e}~\cite[Thm.~1.2.8]{FlOCVo99}
on obtient
$$ 
e(\gothg , A_{I(\eta)}) 
= \lg  (A_{I(\eta)})_{(0)}   
\, e(\gothg , B_{I(\eta)}),
$$  
car  l'id{\'e}al $(0) \subset B_{I(\eta)}$
est le seul id{\'e}al premier minimal  
de $A_{I(\eta)}$. 
On remarque  aussi que 
$ B_{I(\eta)} $ est un  anneau Cohen-Macaulay et que 
$ g_1, \dots, g_n$ 
est un syst{\`e}me de $G$-param{\`e}tres pour $\gothg$, et donc 
\ $  
e(\gothg , B_{I(\eta)})  
= \lg  (B/ \gothg)_{I(\eta)}  = \ell(\eta)$ \ par~\cite[Cor.~1.2.13]{FlOCVo99}. 
Similairement on trouve 
\ $ 
e(\gothg , A_{I(\eta)}) = 
\lg  (A / \gf)_{I(\eta)}   $.

\smallskip 

L'anneau 
$ (A/ \gf)_{I(\eta)} $ est artinien et donc 
\ $
\lg  (A / \gf)_{I(\eta)} = \sum_{\xi \in \psi^{-1} (\eta)} 
\lg  (A / \gf)_{I(\xi)} = \sum_{ \xi \in \psi^{-1} (\eta)}  
\ell(\xi)$.
Finalement 
\ $ \lg  (A_{I(\eta)})_0   = 
[L : K ] = \gamma $; 
en regroupant l'information obtenue on trouve 
$$ 
\sum_{ \xi \in \psi^{-1} (\eta)}  
\ell(\xi) \,=\, \gamma \, \ell(\eta).
$$

\smallskip  

En plus on a  \,$ \widehat{h}_Q(\xi) =  \widehat{h}_P (\eta)$
pour tout $ \xi \in \psi^{-1} (\eta) $, 
car
 \,$  \varphi_\cA (\xi)  = \varphi_\cB (\eta)$.   
On observe aussi que \,$ \NP(g_i)
\subset P$\, et que  \,$ L_{P\cap \Z^n}  = \Z^n$. 
Ainsi on est dans les conditions du cas pr{\'e}c{\'e}dent, d'o{\`u} on d{\'e}duit
$$
\sum_{\xi \in V_0}   \ell(\xi) \, \widehat{h}_Q(\xi)  
\ = \ \gamma \,  
\sum_{\eta \in W_0}   \ell(\eta) \, \widehat{h}_P(\eta)  
\ \le \   
 \gamma \, n!\,\Vol_n(P ) \,  
\sum_{i=1}^n \, h_1 (g_i) 
\  = \  n!\,\Vol_n(Q) \, 
\sum_{i=1}^n \, h_1 (f_i).
$$

\medskip

Le seul cas qui  reste est 
$\dim (Q) < n$. 
Dans cette situation, la vari{\'e}t{\'e} 
torique $ X_\cA \subset \P^N$ est 
aussi de dimension $<n$ et donc les fibres de l'application
$\varphi_\cA: \T^n \to \P^N$ sont toutes de dimension $>0$
gr{\^a}ce au  th{\'e}or{\`e}me de dimension des fibres~\cite[\S~I.6.3]{Sha74}.
Soient \ $\ell_1, \dots , \ell_n \in \Qbarra[y_0, \dots, y_N] $  \ 
les formes lin{\'e}aires correspondant {\`a} $f_1,
\dots, f_n $, de sorte que 
$$
V = \varphi_{\cA}^{-1} \Big( X_\cA \cap 
Z(\ell_1, \dots,  \ell_n) \Big) \ \subset \T^n 
$$
et donc \,$V_0 = \emptyset $. 
Aussi on a 
$\Vol_n(Q) =0$, et donc l'{\'e}nonc{\'e} se r{\'e}duit  
{\`a} l'in{\'e}galit{\'e} triviale $ 0 \le 0$. 
\end{Demo} 

\smallskip

Ce r{\'e}sultat am{\'e}liore 
le cas   
non-mixte du  
th{\'e}or{\`e}me de
Bernstein-Koushnirenko 
arithm{\'e}tique 
d{\^u} {\`a} Maillot~\cite[Cor.~8.2.3]{Maillot00}. 
Avec les notations et les hypoth{\`e}ses du 
th{\'e}or{\`e}me~\ref{2-general}, le r{\'e}sultat de Maillot s'{\'e}crit
\begin{equation} \label{maillot} 
\sum_{\xi \in Z(f_1,\dots, f_n)_0}  \ell(\xi) \, \widehat{h}_Q (\xi) 
\ \le \ 
n!\,\Vol_n (Q) \, \sum_{i=1}^n  \Big( m (f_i) + L(Q)  \Big) ,  
\end{equation} 
o{\`u} \,$m(f_i ) $\, 
d{\'e}signe la mesure de Mahler de $f_i$, 
et \,$L(Q)  $\, est une  constante positive 
associ{\'e}e {\`a} $Q$. 
Le point faible de 
ce r{\'e}sultat
 est son ineffectivit{\'e}, due {\`a}  
la pr{\'e}sence de cette 
constante
$L(Q)$
qu'on ne sait pas
contr{\^o}ler en g{\'e}n{\'e}ral.

\begin{rem} \label{rem-maillot} 

Le seul cas o{\`u} l'on  dispose d'un   
certain contr{\^o}le  de $L(Q)$ 
est quand 
le
polytope $Q$ est {\em absolument simple}; 
dans ce cas
\begin{equation} \label{N(Q)} 
L(Q) \le 3 \, \bigg( \frac{1}{2} \,  \log (n-1) +1  \bigg)
\, N(Q)  
\end{equation}  
pour $n \ge 2$,  o{\`u} $N(Q)$ d{\'e}signe la {\em norme} de
$Q$~\cite[Defn.~8.1.5 et  Prop.~8.1.6]{Maillot00}.
Pourtant, il n'est pas possible de
r{\'e}cup{\'e}rer le th{\'e}or{\`e}me~\ref{2-general} 
{\`a} partir des in{\'e}galit{\'e}s~(\ref{maillot}) et~(\ref{N(Q)})  
m{\^e}me dans cette situation restreinte, 
comme le montre l'exemple suivant: 

\smallskip 

Posons \,$
Q_0 := [0,d]^n  \subset \R^n 
$,  qui est un polytope 
absolument simple de norme 
\,$N(Q_0) = n\, d $. 
Soit \,$f \in \Z[x_1^{\pm 1}, \dots, x_n^{\pm 1}]$\, 
un polyn{\^o}me {\`a} support contenu dans 
$Q_0$. 
Le~\cite[Lem.~1.1]{KrPaSo01} 
implique 
\,$
\log|f| \le  m(f) + n\, d \, \log 2 $,  
en regardant $f$ comme un polyn{\^o}me en $n$ groupes de 1 variable
chacun, de degr{\'e} partiel born{\'e} par $d$ en chaque variable. 
Ainsi 
\begin{eqnarray*} 
h_1(f) & \le & \log||f||_1  \\[2mm] 
& \le & m(f) + n\, d \, \log 2 +  n\, \log(d+1) \\[2mm] 
& < & 
m(f) + 3 \, n\, d  \, \bigg( \frac{1}{2} \,  \log (n-1) +1 \bigg)
= m(f) +  3  \,\bigg( \frac{1}{2} \,  \log (n-1) +1  \bigg)
\, N(Q_0). 
\end{eqnarray*} 
On en d{\'e}duit que le
th{\'e}or{\`e}me~\ref{2-general}  est  (du moins pour cet  exemple)  
 plus
fort que  l'in{\'e}galit{\'e} qui r{\'e}sulte d'appliquer 
l'estimation~(\ref{N(Q)}) dans l'in{\'e}galit{\'e}~(\ref{maillot}).

\end{rem}

\smallskip 

Comme une cons{\'e}quence simple du th{\'e}or{\`e}me~\ref{2-general}, on 
d{\'e}duit l'in{\'e}galit{\'e} de B{\'e}zout
arithm{\'e}tique suivante pour la hauteur de Weil 
des points de
$\T^n$:

\begin{cor} \label{denso}

Soient
\,$f_1, \dots,  f_n \in
\Qbarra[x_1 , \dots, x_n]$\,
des polyn{\^o}mes (ordinaires) 
de degr{\'e} born{\'e} par $d$ 
et  
\  $V:= Z(f_1, \dots, f_n) \subset \T^n$, alors 
$$
\sum_{\xi \in V_0} \ell(\xi) \, \widehat{h}(\xi) 
\ \le \ d^{n-1} \, \sum_{i=1}^n h_1(f_i). 
$$

\end{cor}

\begin{demo} 
Soit $S \subset \R^n$ le simplex standard, 
alors $\Supp(f_i) \subset d \, S$
et donc
$$
\sum_{\xi \in V_0} \ell(\xi) \, \widehat{h}_{d S} (\xi) 
\ \le \ n!\,\Vol_n(d\, S)  \, \sum_{i=1}^n h_1(f_i) 
$$
par le th{\'e}or{\`e}me~\ref{2-general}. 
L'{\'e}nonc{\'e} est impliqu{\'e} par les identit{\'e}s 
\,$\widehat{h}_{d S} (\xi) = d \, \widehat{h} (\xi) $\, (cons{\'e}quence
du lemme~\ref{funtorial}(a))  et 
\,$ \displaystyle \Vol_n(d \, Q ) = \frac{d^n}{n!}$. 
\end{demo}

Il est naturel de se demander si l'on peut obtenir  
une majoration similaire {\`a} celle du
th{\'e}or{\`e}me~\ref{2-general}
pour la 
hauteur $\widehat{h}$ {\`a} la place de  $\widehat{h}_Q$. 
L'exemple suivant montre que la r{\'e}ponse {\`a} cette question  
est n{\'e}gative dans le cas g{\'e}n{\'e}ral:

\begin{exmpl} \label{proj} 
Soit $d, H \in \N^\times$ et posons 
$$
f_1:= x_1- H, \quad  f_2:= x_2 \, x_1^{-d} -H , 
\dots , \quad f_n:= x_n \, x_{n-1}^{-d} -H
\ \ \in
\Z[x_1^{\pm 1}, \dots, x_n^{\pm 1}].
$$
Ceci  est un syst{\`e}me de polyn{\^o}mes de 
Laurent  de hauteur \,$h_1(f_i) = \log (H+1) $\, 
et polytope de Newton 
$$
Q:= \NP(f_1, \dots, f_n) = 
 \Conv \Big( \{0, e_1, 
e_2 - d\, e_1, \dots, e_n - d \, e_{n-1} \} \Big)  \ \subset \R^n 
$$
de volume $1$. 
Pourtant  
$$ 
V:= Z(f_1, \dots , f_n) =\Big\{
\Big(
H, H^{1+ d} , \dots, H^{1+ d+ \cdots +
  d^{n-1}} 
\Big)
 \Big\} \ \subset \T^n
$$ 
est un ensemble {\`a} un seul point 
(en accord avec le th{\'e}or{\`e}me de Koushnirenko g{\'e}om{\'e}trique)
mais de grande hauteur, car 
\,$
\widehat{h}(V) =  \left(1+ d+ \cdots + d^{n-1} \right) \, \log H $.
  
Cependant  
\,$\widehat{h}_Q(V) = \log H $, ce qui est bien en  accord avec la
majoration 
\,$\widehat{h}_Q(V)  \le n \, \log (H+1) $ \ pr{\'e}dite par le
th{\'e}or{\`e}me~\ref{2-general}. 

\end{exmpl}

Toutefois, une telle majoration est valable d{\`e}s que $Q$
contient le simplex standard:

\begin{prop}
Dans les notations et hypoth{\`e}ses du th{\'e}or{\`e}me~\ref{2-general}, 
supposons de plus qu'il 
existe $ b \in \Z^n$ tel que  
$ b+ S \subset Q$, o{\`u} $S$ d{\'e}signe le simplex standard de $\R^n$, 
alors
$$
\sum_{\xi \in Z(f_1,\dots, f_n)_0}  \ell(\xi) \, 
\widehat{h}(\xi) 
\ \le \ n!\,\Vol_n(Q) \,  \sum_{i=1}^n h_1(f_i). 
$$
\end{prop}

\begin{demo}
On a 
\ $ \widehat{h}(\xi) 
= \widehat{h}_{b+ S}(\xi) 
\le  \widehat{h}_{Q}(\xi) $ \ 
6ce au lemme~\ref{funtorial}(b). 
Le reste est une cons{\'e}quence directe  
du th{\'e}or{\`e}me~\ref{2-general}. 
\end{demo}

\medskip

Le r{\'e}sultat suivant g{\'e}n{\'e}ralise le th{\'e}or{\`e}me de Koushnirenko 
arithm{\'e}tique au cas des
intersections arbitraires, 
bien qu'avec la restriction 
$L_\cA = \Z^n$
et une majoration l{\'e}g{\`e}rement plus faible quand on se met 
dans les hypoth{\`e}ses du th{\'e}or{\`e}me~\ref{2-general}. 
On trouve ainsi 
une in{\'e}galit{\'e} de type 
ensembliste, 
dans l'esprit de l'in{\'e}galit{\'e} de B{\'e}zout due {\`a} J.~Heintz~\cite{Heintz83}.

\begin{prop} \label{BK-gral}

Soit  \,$ \cA \subset \Z^n $\, un ensemble fini tel que
 \,$L_\cA = \Z^n$  
et 
\,$f_1, \dots,  f_s \in
\Qbarra[x_1^{\pm 1} , $ $\dots, x_n^{\pm 1}]$\, des polyn{\^o}mes de 
Laurent 
tels que \,$\Supp(f_i) \subset \cA$\,  pour $i=1, \dots, s$, 
alors 
\begin{eqnarray*} 
\deg_\cA \Big(Z(f_1, \dots, f_s) \Big) &  \le &   
\deg_\cA \Big( \div(f_1)  \cdots \div(f_s) \Big) 
\ \le \ 
\Vol (\cA) , \\[2mm]
h_\cA \Big(Z(f_1, \dots, f_s) \Big) 
&  \le  & 
h_\cA \Big(\div(f_1) \cdots  \div(f_s)\Big)  \\[-1mm]
& \le &  
\Vol (\cA) \,
\bigg(  \frac{1}{2} \, (n+1)    
\, \log( \Card( \cA) )
+ 
\sum_{i=1}^s \, h_2 (f_i) 
\bigg). 
\end{eqnarray*}

\end{prop}

\begin{demo}
C'est une cons{\'e}quence directe de 
la relation entre le degr{\'e} et la hauteur du cycle 
\,$ \div(f_1)  \cdots \div(f_s) $\, 
et ceux de son support
\,$  Z(f_1, \dots, f_s) $,  
du lemme~\ref{BK-dot}  et de l'estimation 
pour
$h(X_\cA)$ (corollaire~\ref{altura}). 
\end{demo}

% ----------- Seccion 4 -----------------------------------------------

\typeout{Nullstellensatz} 

\section{Le Nullstellensatz arithm{\'e}tique creux}

\label{Nullstellensatz}

\setcounter{equation}{0}

Finalement nous d{\'e}montrons le th{\'e}or{\`e}me~\ref{nss}. 
Cette d{\'e}monstration s'appuie sur l'in{\'e}galit{\'e}
de type Koushnirenko arithm{\'e}tique suivante,  
pour les intersections dans l'espace affine $\A^n$ 
(comparer avec~\cite[Prop.~2.12]{KrPaSo01}).

\begin{lem} \label{BK-afin}

Soient \,$f_1, \dots, f_s \in \Z[x_1, \dots, x_n]$
pour un certain $n\ge 2$.
Posons $d:= \max_i \deg (f_i) $,  
\ $h:= \max_i h_\mmax(f_i)$ et 
\ $\cA : = 
\Supp \Big(1, x_1, \dots, x_n, f_1, \dots , f_s\Big) \subset \N^n$, 
alors
$$
h \Big( Z(f_1, \dots, f_s) \Big)
\ \le \  
\Vol(\cA) \, \Big(n \, h + 5 \, n\, (n+1) \, \log(d+1) \Big).
$$

\end{lem}

\begin{demo}
On reprend la d{\'e}monstration de~\cite[Prop.~2.12]{KrPaSo01} {\`a} partir de la page 556, ligne 16, 
o{\`u}  l'on remplacera l'application de \cite[Prop.~1.7]{KrPaSo01} par 
le corollaire~\ref{altura}. 
On obtient ainsi
\begin{eqnarray*}
h(V ) &
\le & h(X_\cA) + n\, h\, \deg (X_{\cA}) + 4\, (n+1) \, \log(N+1)
\, \deg (X_{\cA})\\ [1mm]
& \le &  \Vol(\cA) \, 
\Big(n\, h\, + \frac{1}{2} \, (n+1) \, 
\log(N+1) + 4\, (n+1) \, \log(N+1) \Big) \\ [1mm]
& \le & 
 \Vol(\cA) \, 
\Big(n\, h\, + 5\, (n+1) \, \log(N+1) \Big). 
\end{eqnarray*}

Finalement \ $ \log (N+1) \le \log \Big({d+n \choose n} +1 \Big)
\le n \, \log(d+1)$ \ car $n \ge 2$, d'o{\`u} on d{\'e}duit la majoration 
cherch{\'e}e.
\end{demo}

\begin{Demo}{D{\'e}monstration du th{\'e}or{\`e}me~\ref{nss}.--} 

On supposera sans perte de g{\'e}n{\'e}ralit{\'e}  $n,d \ge 2$; les autres cas ont  
{\'e}t{\'e} consid{\'e}r{\'e}s de fa{\c c}on satisfaisante 
dans~\cite[Lem.~3.7 et~3.8]{KrPaSo01}.

\smallskip

D'abord on refait la majoration pour 
$\eta(f_1, \dots,f_s) $ dans~\cite[Lem.~4.9]{KrPaSo01}:
on reprend la notation et la  d{\'e}monstration de ce lemme 
{\`a} partir de la page 590, ligne -2, o{\`u} l'on appliquera le 
lemme~\ref{BK-afin} ci-dessus
{\`a} la place de~\cite[Prop.~2.12]{KrPaSo01}; on en obtient
\begin{eqnarray*}
\eta(f_1, \dots, f_s)  & \le & \Big( n\, \max_i h_\mmax(q_i) 
+ 5 \, n \, (n+1) \, \log (d+1)  \Big) \, \cV\\[1mm] 
& \le & 
\bigg(n\, \Big(h+ \log s + 2\, (n+1)\, \log (d+1)
\Big) + 5 \, n \, (n+1) \, \log (d+1)  \bigg) \, \cV\\[2mm]
& =  & n
\, \cV \, \Big(h+\log s + 7 \, (n+1) \, \log (d+1)  \Big).
\end{eqnarray*}

Pour finir la d{\'e}monstration, on consid{\`e}re la version "intrins{\`e}que{}"
du Nullstellensatz arithm{\'e}tique~\cite[Thm.~2]{KrPaSo01}.
Dans les hypoth{\`e}ses du  th{\'e}or{\`e}me~\ref{nss}, 
on a \,$ \delta \le \cV := \Vol(\cA) $\, 
gr{\^a}ce {\`a}~\cite[Lem.~4.9]{KrPaSo01} et 
$\eta \le  n
\, \cV \, \Big(h+\log s + 7 \, (n+1) \, \log (d+1)  \Big)$
par la majoration ci-dessus. 
On en conclut
\begin{eqnarray*}
h_\mmax(a) , h_\mmax(g_i) & \leq &
(n+1)^2  \, d\, \Big(
2\, \eta +  (h +\log s) \, \delta +
21\,(n+1)^2\, d\, \log (d+1) \, \delta \Big) \\[2mm]
& \le &
(n+1)^2  \, d\, \bigg(
2\, n
\,  \Big(h+\log s + 7 \, (n+1) \, \log (d+1)  \Big) 
\, \cV
+   (h +\log s ) \, \cV   \\[0mm]
&& 
 +
21\,(n+1)^2\, d\, \log (d+1) \, \cV \bigg) \\[2mm]
& \le &
(n+1)^2  \, d\,  \cV \, \Big(
(2\, n +1) \, (h+\log s)
+
28\,(n+1)^2\, d\, \log (d+1) \Big) \\[2mm]
& \le &
2\, (n+1)^3  \, d\,  \cV \, \Big( h+\log s
+
14 \,(n+1)\, d\, \log (d+1) \Big).
\end{eqnarray*}\end{Demo}

% ----------- Referencias -----------------------------------------------

\typeout{Referencias}

\bigskip

\noindent {\sc Mart{\'\i}n Sombra: }
Universit{\'e} de Paris 7, UFR de Math{\'e}matiques,
{\'E}quipe de G{\'e}om{\'e}trie et Dynamique, 
2 place Jussieu, 75251 Paris Cedex 05,
France. \\[1mm] 
{\tt E-mail: sombra@math.jussieu.fr}

\end{document}